\documentclass[hidelinks,onefignum,onetabnum]{siamart251216}

\usepackage{lipsum}
\usepackage{amsfonts}

\usepackage{graphicx}
\usepackage{epstopdf}
\usepackage{subcaption}

\usepackage{algpseudocode}
\usepackage{algorithm}

\usepackage{amssymb}
\usepackage{booktabs}

\usepackage{microtype}

\ifpdf
  \DeclareGraphicsExtensions{.eps,.pdf,.png,.jpg}
\else
  \DeclareGraphicsExtensions{.eps}
\fi

\newsiamremark{remark}{Remark}
\crefname{remark}{Remark}{Remarks}
\AddToHook{env/remark/begin}{\crefalias{theorem}{remark}}
\AddToHook{env/hypothesis/begin}{\crefalias{theorem}{hypothesis}}
\AddToHook{env/example/begin}{\crefalias{theorem}{example}}
\newsiamremark{hypothesis}{Hypothesis}
\crefname{hypothesis}{Hypothesis}{Hypotheses}
\newsiamthm{claim}{Claim}
\newsiamremark{fact}{Fact}
\crefname{fact}{Fact}{Facts}

\theoremstyle{plain}
\theoremheaderfont{\normalfont\sc}
\theorembodyfont{\normalfont\itshape}
\theoremseparator{.}
\theoremsymbol{}
\newtheorem{assumption}{Assumption}

\crefname{assumption}{Assumption}{Assumptions}

\headers{Stochastic semismooth Newton method in Hilbert spaces}{P. Dvurechensky, C. Geiersbach, M. Hintermüller and F. Sauer}

\title{A globally convergent stochastic semismooth Newton method for strongly convex problems in Hilbert spaces 
\thanks{Submitted to the editors November 04, 2026.
\funding{The authors are thankful for support from the Deutsche
Forschungsgemeinschaft (DFG, German Research Foundation) under Germany’s
Excellence Strategy – The Berlin Mathematics Research Center MATH+ (EXC-2046/1, EXC-2046/2, project ID: 390685689). The second and third authors also thank the Deutsche Forschungsgemeinschaft [Projects B02 and B09 in the “Sonderforschungsbereich/Transregio 154
Mathematical Modelling, Simulation and Optimization Using the Example of Gas
Networks”] for support.}}}

\author{
  Pavel Dvurechensky%
  \thanks{Weierstrass Institute for Applied Analysis and Stochastics (WIAS),
  Berlin, Germany
  (\email{pavel.dvurechensky@wias-berlin.de},
   \email{michael.hintermueller@wias-berlin.de},
   \email{felix.sauer@wias-berlin.de}).}%
  \and
  Caroline Geiersbach%
  \thanks{Department of Mathematics, University of Klagenfurt,
  Klagenfurt, Austria
  (\email{caroline.geiersbach@aau.at}).}%
  \and
  Michael Hinterm\"uller\footnotemark[2]
    \thanks{Humboldt-Universit\"at zu Berlin, Berlin, Germany (\email{hint@math.hu-berlin.de}).}
  \and
  Felix Sauer\footnotemark[2]%
}

\usepackage{amsopn}

\ifpdf
\hypersetup{
  pdftitle={A globally convergent stochastic semismooth Newton method for strongly convex problems in Hilbert spaces},
  pdfauthor={P. Dvurechensky, C. Geiersbach, M. Hintermüller and F. Sauer}
}
\fi

\begin{document}

\maketitle

\begin{abstract}
We propose and analyze a globalized stochastic semismooth Newton method for stochastic optimization problems in Hilbert spaces where the objective is strongly convex with Lipschitz continuous, Newton differentiable gradient. Gradient and generalized Hessian information are only available through mini-batching. In each iteration, a full semismooth Newton trial step is accepted or rejected based on a stochastic residual test and rejected steps are replaced by damped steps in the same Newton direction. Throughout, we admit gradient estimators with affinely growing variance, as arises in optimization problems with random partial differential equations as constraints. We prove non-asymptotic global R-linear convergence in expectation as well as almost sure convergence. For sufficiently fast-growing batch sizes, the method transitions almost surely into a regime of pure Newton steps and converges R-superlinearly almost surely. The method is applied to the optimal control of a random elliptic equation featuring a Huber tracking term as well as a random parabolic equation with a Moreau--Yosida penalization of an almost sure state constraint. 
\end{abstract}



\section{Introduction}
Recently, there has been significant progress in modeling, analysis, and numerical methods for optimization problems governed by partial differential equations (PDEs) with uncertain data (see e.g.~\cite{heinkenschloss2025optimization}). In many applications, one seeks a deterministic control $u\in U$ that performs well on average over the realizations of an uncertain parameter $\xi\in \Xi$, which leads to stochastic optimization problems of the form
\begin{equation}\label{eq:intro_problem}
    \min_{u\in U} \ j(u) := \mathbb{E}\big[J(u,\xi)\big],
\end{equation}
where $U$ is a real separable Hilbert space and $J(u,\xi)$ is an uncertain cost functional that measures, for instance, the error of tracking the state of a random PDE driven by the control $u$ towards a desired state. We note that this specific objective $j$ containing the expected value of $J(\cdot, \xi)$ models a risk-neutral scenario. 
Risk-averse formulations have been studied in \cite{alphonse2025risk, kouri2016risk, kouri2018existence}. We are interested in a class of problems \cref{eq:intro_problem} exhibiting two structural features that are characteristic of PDE-constrained optimization under uncertainty. 
First, the expectation is not available in closed form. Thus, every evaluation of $j$ or its derivatives must be approximated by sampling, and each sample requires computing the solution of one or more PDEs. Second, the problem is nonsmooth in the sense that the objective $j$ is continuously differentiable, but the gradient $g = \nabla j$ is Newton (slantly) differentiable. This renders $j$ a so-called $SC^1$ function \cite{vonheusinger2008SC1OptimizationReformulations} and second-order information is only available through generalized differentiation. For instance, Moreau--Yosida penalizations of state constraints or Huber-type tracking terms lead to this structure.

A numerical treatment of \eqref{eq:intro_problem} requires a discretization of the PDE as well as the sample space $\Xi$. Frequently used sampling strategies include sample average approximation (SAA) and (quasi-)Monte Carlo methods on the one hand and stochastic approximation on the other (compare \cite{heinkenschloss2025optimization} and the references therein). This work focuses on stochastic approximation, which was first investigated for PDE-constrained problems in Hilbert spaces in \cite{geiersbach2019ProjectedStochasticGradients,geiersbach2020stochastic,Martin2021,Martin2021a} and for inverse problems in \cite{Dippon2023}. Accelerated variants are found in \cite{Toraman2023}. 
Stochastic first-order methods are robust and easy to implement, but they converge only sublinearly \cite{bottou2018optimization} and their performance deteriorates with the conditioning of the problem. This is particularly severe for path-following approaches in constrained optimization, where a penalty parameter is driven to infinity to enforce feasibility \cite{hintermuller2009MoreauYosidaRegularization, hintermuller2006pathfollowing}. In the deterministic setting, semismooth Newton methods \cite{hintermuller2010SemismoothNewtonMethods,hintermuller2002PrimalDualActiveSet,kunisch2005GeneralizedNewtonMethods,ulbrich2011semismooth} have been used successfully to solve $SC^1$ problems. They use generalized derivatives of the nonsmooth optimality system for Q-superlinear local convergence in function spaces.

Newton-type methods with sampled derivatives have recently received considerable attention in the finite-dimensional setting. For smooth problems, subsampled Newton methods were analyzed in \cite{bollapragada2019ExactInexactSubsampled}, where global R-linear and local R-superlinear convergence in expectation was established. Convergence was established under a uniformly bounded gradient variance \cite[A2]{bollapragada2019ExactInexactSubsampled} and the local rate is obtained for initial points sufficiently close to the solution \cite[Theorem 2.6]{bollapragada2019ExactInexactSubsampled}. However, the analysis is finite-dimensional and requires twice continuously differentiable objectives as well as uniformly bounded gradient variance. The former excludes the nonsmooth structures described above; the latter fails whenever the random parameter is not additively separable from the control in the PDE, see \cref{sec:SSN}.

For nonsmooth problems, a stochastic semismooth Newton method was introduced by \cite{milzarek2019stochastic} for nonconvex, composite objectives.
There, a semismooth Newton trial step is accepted under growth conditions with a stochastic proximal gradient step as a fallback; the expected residual is shown to converge to zero, but without an associated convergence rate. The accompanying local theory in \cite{milzarek2022local} establishes that the method eventually only takes pure Newton steps and converges R-linearly or R-superlinearly with high probability, provided the solution is an accumulation point of the generated sequence \cite[Assumption 4.2]{milzarek2022local}. However, we note that the proofs rely on matrix concentration estimates \cite[Lemma 3.4]{milzarek2022local} with constants depending on the dimension $d$ which are inherited by the generalized Hessian batch sizes \cite[Theorem 4.5]{milzarek2022local}. Such estimations cannot simply be carried over to the infinite dimensional setting since even for a discretized PDE, $d$ grows under mesh refinement. 
In the deterministic Hilbert space setting, a globalized inexact semismooth Newton method for strongly convex problems, applied to the continuously differentiable Fenchel dual, was recently proposed in \cite{wachsmuth2026globalized}.

To the best of our knowledge, Newton-type methods with both sampled gradients and sampled Hessians have so far only been analyzed in finite dimensions with fast local convergence depending on the initial point or the generated sequence \cite{bollapragada2019ExactInexactSubsampled,milzarek2022local}. We close this gap in Hilbert spaces for sampled gradients and generalized Hessians, where the gradient variance is allowed to grow affinely. We prove a nonasymptotic global R-linear rate in expectation from any initial point, with explicit contraction factor and batch-size thresholds, already for generalized Hessian batches of size one. For the semismooth method of \cite{milzarek2019stochastic}, no global rate is available. Moreover, the transition to full Newton steps and R-superlinear convergence follows almost surely from the global phase instead of being assumed. Batch sizes and rates depend only on problem constants.

\hspace{-0.5em}We propose and analyze a globalized stochastic semismooth Newton (SSN) method for problems \eqref{eq:intro_problem} with $\mu$-strongly convex objective and $L$-Lipschitz gradient in a general, possibly infinite-dimensional Hilbert space. In each iteration, we set 
\[u_{n+1} = u_n + \tau_n d_n \quad \text{with $d_n$ solving} \quad H_{S_n}(u_n)[d_n] = - G_{X_n}(u_n).\]
The operators $G_{X_n}$ and $H_{S_n}$ are mini-batched estimators of the gradient and generalized Hessian, respectively. 
A stochastic acceptance test, closely related to the growth conditions of \cite{milzarek2019stochastic}, determines the step size $\tau_n>0$ by comparing the residual at the full-step trial point ($\tau_n = 1$), evaluated on a fresh test batch, to the residual at the current iterate. If the test passes, the full-step is accepted; otherwise, the method falls back to a damped step. 

In \cref{sec:SSN}, we outline the method and state preliminary assumptions.  
Then, in \cref{sec:convergence_analysis}, we prove nonasymptotic global R-linear convergence in expectation with an explicit contraction factor and explicit lower thresholds for gradient batch sizes and Hessian batches of size one (\cref{thm:global_linear_combined}). We show almost sure convergence of the iterates (\cref{lem:almost_sure_convergence}), and the almost sure transition to full Newton steps with R-superlinear convergence
(\cref{thm:local_r_superlinear}). Entry into the local regime of fast convergence is thus a consequence of the global phase instead of an assumption on the initialization as in \cite{bollapragada2019ExactInexactSubsampled} or an assumption on the generated sequence \cite{milzarek2022local}.

In \cref{sec:numerics}, we illustrate the method on two PDE-constrained model problems. The first test example involves a random linear elliptic state equation with a Huber tracking objective. Here, globalization is demonstrated. Initialized outside of the quadratic region of the Huber loss, the full-step method does not converge, while the globalized method takes a few damped steps and then switches to full Newton steps to quickly reach the sampling-induced noise floor. The second test example involves a random linear parabolic equation with a Moreau--Yosida penalization of an almost sure state constraint. In both cases, we compare against stochastic gradient descent, both per iteration and in terms of number of PDE solves.

\section{Stochastic semismooth Newton method}\label{sec:SSN}

\subsection{Preliminaries}
\paragraph{Spaces and operators}
Let $U$ be a real, separable Hilbert space with inner product $\langle\cdot,\cdot\rangle$ and induced norm $\|\cdot\|$. We denote by $\mathcal{L}(U)$ the space of bounded linear operators from $U$ to $ U$ with operator norm $\|A\|:= \sup_{\|u\|=1}\|A[u]\|.$ We call $A \in \mathcal{L}(U)$ self-adjoint if $\langle A[v], w\rangle = \langle  v, A[w]\rangle$ for all $v,w \in U$, and coercive if there exists $c>0$ such that $\langle A[v], v\rangle \geq c \|v\|^2$ for all $v \in U.$ 
Lastly, we use the convention $\mathbb{N}_0 := \mathbb{N} \cup \{0\}$ and $\mathbb{R}^+ := [0, \infty).$ 

\paragraph{Differentiability}
For a function $j\colon U \to \mathbb R$ and $u \in U$, we denote by $Dj(u) \in \mathcal{L}(U, \mathbb R)=: U^\ast$ its Fréchet derivative. The gradient $\nabla j$ is its Riesz representation, i.e.,~$\langle \nabla j(u), v \rangle = \langle D j(u), v\rangle_{U^\ast, U}$ for all $v \in U.$ We denote by $C_L^{1,1}(U)$ the class of continuously Fréchet differentiable functions mapping $U$ to  $\mathbb{R}$ with $L$-Lipschitz gradient, i.e.,~there exists $L>0$ such that for all $u,v \in U$
\begin{equation}\label{eq:def_lipschitz}
    \|\nabla j(u) - \nabla j(v) \| \leq L \|u-v\|.
\end{equation} 
Functions in $C_L^{1,1}(U)$ satisfy the descent lemma (cf.~\cite[Theorem 18.15]{bauschke2011ConvexAnalysisMonotone}); specifically, for all $u, v \in U$ we have
\begin{equation}\label{eq:def_descent_lemma}
    j(v) \leq j(u) + \langle \nabla j(u), v-u \rangle + \frac{L}{2}\|v-u\|^2.
\end{equation}
Further, we call $j$ an $SC^1$ function if it is continuously Fréchet differentiable with Newton differentiable Fréchet derivative (cf. \cite{vonheusinger2008SC1OptimizationReformulations}).

\paragraph{Probability}
Throughout, $(\Omega, \mathcal{F}, \mathbb{P})$ will denote a complete probability space and $\Xi$ is a complete, separable metric space, such that $\xi\colon \Omega \to \Xi$ is measurable. We often suppress the dependence on $\omega \in \Omega$ and use the abusive notation $\xi \in \Xi$ for a realization. For a set $B$ of samples, we denote its cardinality by $|B|$.

Since $U$ is separable, weak and strong measurability of maps from $\Omega$ to $U$ coincide and we speak of $U$-valued random variables. We denote by $\mathbb E[f] := \int_\Omega f(\omega)\, \mathrm{d} \mathbb{P}(\omega) \in U$  
the expectation of a $U$-valued random variable $f$ in the Bochner sense. Expectations of operator-valued maps are understood pointwise: for
$h \colon \Xi \to \mathcal{L}(U)$ we let $\mathbb{E}[h(\xi)]$ denote the
operator defined by $\mathbb{E}[h(\xi)][d] := \mathbb{E}\big[ h(\xi)[d] \big]$ for all $d \in U.$

A filtration is an increasing sequence $\{\mathcal{F}_n\}$ of $\sigma$-algebras of $\mathcal{F}$, i.e.,~$\mathcal{F}_1 \subset \mathcal{F}_2 \cdots \subset \mathcal{F}.$  A $U$-valued random variable $v$ is called $\mathcal{F}_n$-measurable if $v^{-1}(B) \in \mathcal{F}_n$ for every Borel set $B \subseteq U.$ 
For a collection of $U$-valued random variables $\{v_n \colon \Omega \to U : n \in \mathbb{N}\}$, the natural filtration is the filtration generated by the sequence $\{v_n\}$ and is defined by $\mathcal{F}_n= \sigma(\{v_1,...,v_n\})$. For a $U$-valued random variable $v$, its conditional expectation $\mathbb{E}[v \mid \mathcal{F}_n]$ is an $\mathcal{F}_n$-measurable random variable, which by definition satisfies $\int_A \mathbb{E}[v \mid \mathcal{F}_n] \, \mathrm{d}\mathbb{P}(\omega) = \int_A v(\omega) \, \mathrm{d}\mathbb{P}(\omega)$
for every $A \in \mathcal{F}_n$. For an event $A \in \mathcal{F}$, we write $\mathbf{1}_A$ for
its indicator function, $A^c$ for its complement, and
$\mathbb{P}(A \mid \mathcal{F}_n) := \mathbb{E}[\mathbf{1}_A \mid \mathcal{F}_n]$ for the conditional probability.

\subsection{Problem statement, assumptions and algorithm}
In the following, we consider the stochastic optimization problem \cref{eq:intro_problem}
for a continuously differentiable $j \colon U \rightarrow \R$.
Let $g\colon U \to U$ be the gradient operator 
\begin{equation*}\label{eq:def_g}
g(u) :=  \nabla j(u) = \mathbb E[\nabla_u J(u,\xi)].
\end{equation*} 
We introduce the single sample gradient $G\colon U \times \Xi\to U, \ (u, \xi) \mapsto G_\xi(u) := \nabla_u J(u, \xi) $ such that $g(u) = \mathbb{E}[G_\xi(u)]$ for all $u \in U$. Throughout, capital letters denote sampled quantities and the subscript indicates the sample or batch used. Lowercase letters denote the corresponding expected quantities. 

In the following, we frequently invoke the following assumptions. 
\begin{assumption}[Strongly convex objective with Lipschitz gradient]\label{ass:lipschitz_mu_convex}
    We assume $j \in C_L^{1,1}(U)$, i.e.,~the gradient $g = \nabla j$ satisfies \eqref{eq:def_lipschitz}. 
    We further assume that $j$ is $\mu$-strongly convex, i.e.,~there exists $\mu > 0$ such that for all $u,v \in U$ we have
    \begin{equation}\label{eq:def_strongly_convex}
      j(v) \geq j(u) + \langle g(u) , v-u\rangle + \frac{\mu}{2}\|u-v\|^2.  
    \end{equation}
\end{assumption}

We suppose that $u^\ast\in U$ is a solution of $g(u)=0$, i.e.,~$u^\ast$ solves the minimization problem \eqref{eq:intro_problem}. 
Under \cref{ass:lipschitz_mu_convex}, $u^\ast$ is unique.

To formulate the SSN method, we employ the following notion of
Newton (slant) differentiability at $u^\ast$; cf.~\cite{hintermuller2010SemismoothNewtonMethods,hintermuller2002PrimalDualActiveSet} or \cite{ulbrich2011semismooth} for a related notion.  

\begin{assumption}[Newton differentiability of the gradient]\label{ass:newton_differentiability}
We assume that there exists a mapping $h\colon U \to  \mathcal{L}(U)$
such that at $u^\ast \in U$,
\begin{equation}\label{eq:slanting_function}
\lim_{\|s\|\to 0}
\frac{\|g(u^\ast+s)-g(u^\ast)-h(u^\ast+s)[s]\|}{\|s\|}=0.
\end{equation}
We call $g$ Newton (or slantly) differentiable at $u^\ast$, and $h$ is referred to as a generalized derivative of $g$. 
\end{assumption}
Since $g=\nabla j$, we will also refer to $h(u)$ as a generalized Hessian of $j(u)$ and assume that there exists an estimator $H\colon U \times \Xi \to \mathcal{L}(U),\ (u, \zeta) \mapsto H_\zeta(u)$ with $h(u) = \EE{H_\zeta(u)}$ for all $u \in U.$ We further assume that $(u, \xi) \mapsto G_\xi(u)$ and $(u,d,\zeta)\mapsto H_\zeta(u)[d]$
are Borel measurable.

We will require estimators of the gradient and (generalized) Hessian built from sample batches. At iteration $n$, the algorithm draws 
\[X_n := \{\xi_1^n,...,\xi^n_{|X_n|}\}, \quad S_n := \{\zeta^n_1,...,\zeta^n_{|S_n|}\},\quad \tilde{X}_n := \{\tilde \xi^n_1,...,\tilde\xi^n_{|\tilde X_n|}\},\] 
and forms the mini-batch estimators 
\begin{equation*}\label{eq:canonical_estimators}
    G_{X_n}(u) = \frac{1}{|X_n|}\sum_{\xi \in X_n} G_\xi(u), \qquad H_{S_n}(u) = \frac{1}{|S_n|}\sum_{\zeta \in S_n} H_\zeta(u),
\end{equation*}
and analogously $G_{\tilde{X}_n}$. Throughout, we assume that all samples are independent and identically distributed (i.i.d.). 
We define the filtrations
\begin{equation}\label{eq:filtrations}
    \mathcal{F}_n:=\sigma(X_k, S_k, \tilde{X}_k : k \leq n-1), \qquad \mathcal{F}_n^+:= \sigma(X_{k}, S_{k}, \tilde{X}_{k-1} : k \leq n),
\end{equation} 
so that $\mathcal F_n \subseteq \mathcal F_n^+ \subseteq \mathcal F_{n+1}$; this implies that $X_n, S_n$ are independent of $\mathcal{F}_n$ and $\tilde{X}_n$ is independent of $\mathcal{F}_n^+.$ 

\begin{assumption}[Gradient estimator with affinely growing variance]\label{ass:gradient_variance}
    We assume that there exist constants $\nu_0, \nu_1 \geq 0$, such that for all $u \in U$,    
\begin{equation}\label{eq:global_grad_variance}
\mathbb{E}\left[\|G_\xi(u) - g(u)\|^2 \right] \leq {\nu_0^2 + \nu_1^2 \|u- u^\ast\|^2}.
    \end{equation}
\end{assumption}

\begin{assumption}
    [Generalized Hessian estimator with bounded variance]\label{ass:hessian_variance}
There exists a constant $\sigma \geq 0$  such that for all $u, d \in U$,
\begin{equation}\label{eq:global_hess_variance}
          \mathbb{E}\left[\|(H_\zeta(u) - h(u))[d]\|^2 \right] \leq {\sigma^2 \|d\|^2}.
    \end{equation}
\end{assumption}

We note that the bound in \cref{ass:gradient_variance} is similar to a second moment gradient bound as in \cite[Assumption 3.4 (iii)]{geiersbach2019ProjectedStochasticGradients} or the variance condition in \cite[Assumption 4.3]{bottou2018optimization}.
Further, we note that \cref{ass:gradient_variance} accounts for applications in PDE-constrained optimization where randomness entering the PDE is not linearly separable from the control. Thus, the affinely growing variance \eqref{eq:global_grad_variance} sets us apart from related second-order methods such as \cite{bollapragada2019ExactInexactSubsampled,milzarek2019stochastic,milzarek2022local}, which assume uniformly bounded variance. Further, we center \cref{ass:gradient_variance} around $u^\ast$, since only distances to $u^\ast$ enter the analysis. Because $\|u\|^2 \leq 2\|u-u^\ast\|^2 + 2\|u^\ast\|^2$, the bound \eqref{eq:global_grad_variance} is equivalent (up to the constants $\nu_0, \nu_1$) to a bound of the form $\nu_0^2 + \nu_1^2\|u\|^2$.

\begin{assumption}[Uniform bounds on the sampled generalized Hessian]\label{ass:hessian_bounds}
We assume there exist constants $0 < m \leq M < \infty$ such that for every $ u \in U$, and almost every $\zeta \in \Xi,$ the operator $H_{\zeta}(u)$ is self-adjoint and satisfies 
$m \|v\|^2  \leq  \langle H_\zeta(u)[v], v \rangle \leq M\|v\|^2$ for all  $v \in U$. For $n \in \mathbb{N}_0$, this implies self-adjointness of the batched operator $H_{S_n}(u)$ and 
\begin{equation}\label{eq:hessian_bounds}
      m \|v\|^2  \leq  \langle H_{S_n}(u)[v], v \rangle \leq M\|v\|^2 \qquad \forall v \in U.
    \end{equation}
\end{assumption}

Throughout the paper we set 
\begin{equation}\label{eq:def_rho_theta}
        \rho := 1 - \frac{m\mu}{LM}, \qquad \theta := \sqrt{\frac{\rho\mu}{L}}.
\end{equation}
Since $\tfrac{m}{M} \leq 1$ and $\tfrac{\mu}{L} \leq 1$, we have $\rho \in [0,1)$, and $\theta^2 = \rho\tfrac{\mu}{L} \leq \rho$ gives $\theta \in [0,1)$. Moreover, $\rho = 0$ holds if and only if $m = M$ and $\mu = L$. We exclude this case as a standing assumption, so that $\rho, \theta \in (0,1)$.

We can now state the algorithm. 
\begin{algorithm}[ht]
\caption{Globalized Stochastic Semismooth Newton Method}
\label{alg:globalized_SSN}
\begin{algorithmic}[1]
\Require initial iterate $u_0\in U$, problem constants $m,M, \mu, L$ with $\theta$ as defined in \eqref{eq:def_rho_theta}, batch sizes $|X_n|, |S_n| \geq 1$ and a slack sequence $\delta_n\geq 0$ for $n \in \mathbb{N}_0$
\For{$n=0,1,2,\dots$}
  \State draw i.i.d.~batches of samples $X_n,S_n, \tilde{X}_n$ with $|X_n| = |\tilde{X}_n|$ 
  \State solve $\;H_{S_n}(u_n)\,[d_n]=-\,G_{X_n}(u_n)\;$ for $d_n\in U$
  \State set the full-step trial point $\;z_n := u_n + d_n$
  \State evaluate the stochastic test
  \begin{equation}\label{eq:ast}\tag{\mbox{$\ast$}}
       \big\|G_{\tilde X_n}(z_n)\big\|
       \;\le\;\theta\,\big\|G_{
       { X_n}}(u_n)\big\|+\delta_n
  \end{equation}
    
  \If{$(\ast)$ holds}
    \State $u_{n+1} := z_n$ 
  \Else
    \State $u_{n+1} := w_{n} := u_n+\frac{m}{L}\,d_n$

  \EndIf
\EndFor
\end{algorithmic}
\end{algorithm}
\begin{remark}
\Cref{alg:globalized_SSN} requires the constants $m, M, \mu, L$, which enter the test parameter $\theta$ and the damped step size $m/L$. However, only bounds are needed, so underestimating $m, \mu$ or overestimating $M, L$ preserves the convergence results at the price of a stricter test \eqref{eq:ast} and shorter damped steps.
\end{remark}

\subsection{Basic properties}
\cref{ass:lipschitz_mu_convex} implies the following inequalities. 
\begin{lemma}
 Under \cref{ass:lipschitz_mu_convex}, for any $u \in U$, we have
 \begin{align}
 \label{eq:S1}
     \frac{1}{2L}\|g(u)\|^2 \leq j(u) - j(u^\ast ) \leq \frac{1}{2\mu}\|g(u)\|^2, \\
\label{eq:S2}
    \mu \|u-u^\ast\| \leq \|g(u)\| \leq L \|u-u^\ast\|,  
\\
\label{eq:S3}
    \frac{\mu}{2}\|u-u^\ast\|^2 \leq j(u)- j(u^\ast). 
 \end{align}

\end{lemma}
\begin{proof}The proof is the same as in the finite-dimensional setting (see e.g. \cite[Sections 2.1.1, 2.1.3]{nesterov2018lectures}).
 \end{proof}

\begin{lemma}[Properties of the generalized Hessian]
    Let $H_{S_n}$ satisfy \cref{ass:hessian_bounds}. Then $H_{S_n}$ is boundedly invertible with 
\begin{equation}
    \|H_{S_n}^{-1}\|_{\mathcal{L}(U)} \leq \frac{1}{m}  \quad \text{and} \quad 
        \frac{1}{M}\|v\|^2 \leq \langle v, H_{S_n}^{-1}[v]\rangle, \quad \forall v \in U.  \label{eq:H_inv_bounds} 
    \end{equation}
    Further, with $\tau:= \frac{m}{L}$, we have the bound
\begin{equation}\label{eq:descent_lemma_hessian_bound}
        - \tau \langle v, H_{S_n}^{-1}[v]\rangle + \frac{L \tau^2}{2}\|H_{S_n}^{-1}[v]\|^2 \leq - \frac{\tau}{2M}\|v\|^2, \quad \forall v \in U.
    \end{equation}
\end{lemma}
\begin{proof} 
    By \cref{ass:hessian_bounds}, $H_{S_n}$ is self-adjoint with spectrum in $[m,M]$ . Hence, $H_{S_n}$ is boundedly invertible with spectrum $[\frac{1}{M}, \frac{1}{m}]$ which yields \eqref{eq:H_inv_bounds} (see also \cite[Chapter 6]{brezis2011FunctionalAnalysisSobolev}).
    Since $\|H_{S_n}^{-1}[v]\|^2 \leq \frac{1}{m}\langle v, H_{S_n}^{-1}[v]\rangle$ and $\tau = \frac{m}{L}$, the left-hand side of \eqref{eq:descent_lemma_hessian_bound} is bounded by $-\frac{\tau}{2}\langle v, H_{S_n}^{-1}[v]\rangle \leq -\frac{\tau}{2M}\|v\|^2$.
\end{proof}

\begin{remark}[Variance of mini-batch estimators]\label{lem:variance-mini-batch}
    Let \cref{ass:gradient_variance} and \cref{ass:hessian_variance} hold. Then for every $n \in \mathbb{N}_0$ and $u,d \in U$, we have
    \begin{equation*}
         \mathbb{E}\left[\|G_{X_n}(u) - g(u)\|^2 \right] \leq \frac{\nu_0^2 + \nu_1^2 \|u- u^\ast\|^2}{|X_n|}, \quad  \mathbb{E}\left[\|(H_{S_n}(u) - h(u))[d]\|^2 \right] \leq \frac{\sigma^2}{|S_n|}\|d\|^2,
    \end{equation*}
  and the same bound holds for $G_{\tilde{X}_n}$.
\end{remark}

\begin{remark}\label{rem:conditional-variance-minibatch}
We fix $n \in \mathbb{N}$ and let $u,d\colon \Omega \to U$ be $\mathcal{F}_n$ measurable. Further, the maps $(u, X_n) \mapsto \|G_{X_n}(u) - g(u)\|^2$ and $(u, d , S_n)\mapsto \|(H_{S_n}(u) - h(u))[d]\|^2$ are measurable. As $X_n$ and $S_n$ are independent of $\mathcal{F}_n$, we can use \cite[Proposition 1.12]{da2014stochastic} to fix $u, d$ in \cref{lem:variance-mini-batch} which yields
\vspace{-1em}
\begin{align}
    \label{eq:cond_grad_var}
        \mathbb{E}\left[\|G_{X_n}(u) - g(u)\|^2 \mid \mathcal{F}_n \right] &\leq \frac{\nu_0^2 + \nu_1^2 \|u-  u^\ast\|^2}{|X_n|},\\
        \label{eq:cond_hess_var}
        \mathbb{E}\left[\|(H_{S_n}(u) - h(u))[d]\|^2 \mid \mathcal{F}_n \right] &\leq \frac{\sigma^2}{|S_n|}\|d\|^2,
\end{align}
    and the same bound holds for $G_{\tilde{X}_n}$, conditioned on $\mathcal{F}_n^+$ for every $\mathcal{F}_n^+$-measurable $u$. 
\end{remark}

The coercivity constant $m$ as well as the notion of semismoothness in \eqref{eq:slanting_function} lead to the following properties. 

\begin{remark}[Basic properties of the semismoothness]\label{rem:basic_SS}
    For $r>0$, we set %
\begin{equation}\label{eq:semismooth_modulus}
    \omega(r):= \sup_{0 < \|s\| \leq r} \frac{\|g(u^\ast+s)- g(u^\ast) - h(u^\ast+s)[s]\|}{m\|s\|}.
\end{equation}
     Then $\omega$ is 
    \begin{enumerate}[(i)]
        \item monotonically non-decreasing with $\omega(r_1) \leq \omega(r_2)$ for $0 < r_1 \leq r_2$,
        \item vanishing at the origin with $\lim_{r \to 0} \ \omega(r) =0$,
        \item finite with $\omega(r) \leq (L+M)/m$ for all $r>0.$
    \end{enumerate}
\end{remark}
We refer to $\omega$ as the semismoothness modulus at $u^\ast$ and note that ``semismoothness of order $\alpha$ at $u^\ast$" as stated in \cite[Definition 3.1]{ulbrich2011semismooth} implies $\omega(r) = \mathcal{O}(r^\alpha)$ as $r \to 0$.

We close this section with some estimates on the noise level of the iterates. Throughout the following, we note that $u_n$ is $\mathcal{F}_n$-measurable and $z_n, w_n$ are $\mathcal{F}_n^+$-measurable (see, e.g., \cite[Lemma 2.5]{geiersbach2023StochasticGradientMethod}).

\begin{remark}[Iterate and gradient noise bounds]\label{rem:iterate_noise_bound} 
Define the $\mathcal{F}_n$-measurable quantity
\vspace{-1em}
 \begin{equation}\label{eq:def_Vn}
    V_n := \nu_0^2 + \frac{2\nu_1^2}{\mu}\left(j(u_n) - j(u^\ast)\right).
\end{equation}
By \eqref{eq:S3} we have $\nu_0^2 + \nu_1^2\|u_n - u^\ast\|^2 \leq V_n$, and hence
 by \cref{rem:conditional-variance-minibatch},
 \begin{equation}\label{eq:variance_trick}
     \mathbb{E}\big[\|G_{X_n}(u_n)-g(u_n)\|^2 \mid \mathcal{F}_n\big]\leq
     \frac{\nu_0^2 + \nu_1^2 \|u_n- u^\ast\|^2}{|X_n|} \leq \frac{V_n}{|X_n|}.
 \end{equation}
\end{remark}

Similar to \cref{rem:iterate_noise_bound}, we first formulate a bound on the noise at the full-step trial point $z_n$ in terms of the noise level $V_n$ at $u_n$.
\begin{lemma}[Trial-point noise level]\label{lem:initial_zn_bound}
Let \cref{ass:lipschitz_mu_convex,ass:gradient_variance,ass:hessian_bounds} hold. For the full-step trial point $z_n = u_n - H_{S_n}^{-1}(u_n)[G_{X_n}(u_n)]$, we have
\vspace{-1em}
    \begin{equation}\label{eq:initial_zn_bound}
        \nu_0^2 + \nu_1^2\, \mathbb{E}[\|z_n- u^\ast\|^2 \mid \mathcal{F}_n] \leq C_z V_n,
        \qquad
        C_z := 2 + \frac{4(\mu L + \nu_1^2)}{m^2}.
    \end{equation}
\end{lemma}
\begin{proof}
    By the triangle inequality and \eqref{eq:H_inv_bounds}, it holds that
    \vspace{-1em}
   \begin{align*}
  \|z_n - u^\ast\|^2
  &\le 2\|u_n - u^\ast\|^2
     + \frac{2}{m^2}\|G_{X_n}(u_n)\|^2 \\
  &\le 2\|u_n - u^\ast\|^2
     + \frac{4}{m^2}\Big(\|g(u_n)\|^2
     + \|G_{X_n}(u_n) - g(u_n)\|^2\Big).
\end{align*}
    We bound each term against $V_n$. By \cref{rem:iterate_noise_bound}, we have 
    $ \nu_0^2 + 2\nu_1^2\|u_n- u^\ast\|^2 \leq 2 V_n .$ By \cref{eq:S1} and \eqref{eq:def_Vn}, $\nu_1^2\|g(u_n)\|^2 \leq 2L\nu_1^2 (j(u_n)-j(u^\ast)) \leq \mu L V_n.$
    Lastly, by \eqref{eq:variance_trick} and $|X_n| \geq 1$, we have $\nu_1^2 \mathbb{E}[\|G_{X_n}(u_n) - g(u_n)\|^2 \mid \mathcal{F}_n] \leq \frac{\nu_1^2 V_n}{|X_n|} \leq \nu_1^2 V_n.$
    Since $u_n$ is $\mathcal{F}_n$-measurable, we can collect the three bounds and obtain 
    \begin{align*}
        \nu_0^2 + \nu_1^2 \mathbb{E}[\|z_n- u^\ast\|^2 \mid \mathcal{F}_n]
    &\leq \Big(2 + \frac{4\mu L}{m^2}
        + \frac{4\nu_1^2}{m^2}\Big) V_n
  = C_z V_n.
    \end{align*}

\end{proof}

\section{Convergence analysis}\label{sec:convergence_analysis}

\subsection{Global linear convergence with adaptive steps}
We now move to the convergence analysis of \cref{alg:globalized_SSN}.

\begin{proposition}[Global objective decrease]\label{prop:global_objective_decrease}
    Let \cref{ass:lipschitz_mu_convex,ass:gradient_variance,ass:hessian_bounds} hold. Let the iterates $\{u_n\}$ be generated by \cref{alg:globalized_SSN}. 
    Then,
    \begin{equation}\label{eq:global_linear}
      \mathbb{E}[j(u_{n+1})-j(u^\ast) \mid \mathcal{F}_n]\leq \rho_n (j(u_n)-j(u^\ast))+e_n ,
   \end{equation}
   \vspace{-2em}
   \begin{align*}
      \text{where }\;\;\rho_n := \sqrt{\rho}
        + \frac{2\nu_1^2 K}{\mu\,|X_n|},
      \qquad
      e_n := \frac{\nu_0^2 K}{|X_n|}
        + \frac{3\,\delta_n^2}{2\mu(1-\sqrt{\rho})},
   \end{align*}
   with constant $ K:= \frac{1}{2(1-\sqrt{\rho})}\left( \frac{3 (\theta^2 + C_z)}{\mu} + \frac{1}{L}\right).
  $
   Moreover, $\rho_n\leq \bar{\rho}:=\frac{1+\sqrt{\rho}}{2}<1$ with
   \begin{equation}\label{eq:X_min}
      |X_n| \geq X_{\min} := \frac{4\nu_1^2 K}{\mu(1-\sqrt{\rho})}.
   \end{equation}
    
\end{proposition}

\begin{proof}
We introduce the notation 
\begin{align*}
    \Delta_n := j(u_n) - j(u^\ast), \quad
    \pi_n := G_{X_n}(u_n) - g(u_n), \quad 
    \tilde{\pi}_n := G_{\tilde{X}_n}(z_n) - g(z_n),
\end{align*}
and define the event $A_{n} := \{  \eqref{eq:ast} \text{ holds}\}.$
Since $u_{n+1} = z_n$ on $A_n$ and $u_{n+1} = w_n$ on $A_n^c$, we can write   \begin{equation}\label{eq:split}
      \Delta_{n+1}
      =\mathbf 1_{A_n}\big(j(z_n)-j(u^\ast)\big)+\mathbf 1_{A_n^c}\big(j(w_n)-j(u^\ast)\big),
   \end{equation}
   and we bound the two terms separately. Recall the definitions $\rho = 1- \frac{m \mu}{LM}$ and $\theta = \sqrt{\tfrac{\rho \mu}{L}}$ from \eqref{eq:def_rho_theta}. Throughout the proof, we use Young's inequality to estimate $(a+b)^2 \leq (1+\zeta)a^2 + (1+\zeta^{-1})b^2$ given $a,b\geq 0$ and $\zeta>0$. %
In particular, with $\zeta :=\frac{1}{\sqrt{\rho}}-1>0$, it follows that
\begin{equation}\label{eq:young_2}
    (a+b)^2 \leq \frac{1}{\sqrt{\rho}}a^2 + \frac{1}{1-\sqrt{\rho}}b^2.
\end{equation}
    
    \medskip\noindent\emph{Part 1: full-step $z_n$.} 
    On $A_n$, the test \eqref{eq:ast} and the triangle inequality give 
    \begin{align*}
        \|g(z_n)\|\leq\|G_{\tilde{X}_n}(z_n)\|+\|\tilde{\pi}_n\|
       &\leq\theta\|G_{{X}_n}(u_n)\|+\delta_n+\|\tilde{\pi}_n\| \leq \theta\|g(u_n)\|+R_n
    \end{align*}
    with $R_n:= \delta_n+\theta \|\pi_n\|+\|\tilde{\pi}_n\|.$
   Squaring both sides in order to employ \eqref{eq:S1} and using \eqref{eq:young_2} with $a=\theta\|g(u_n)\|$, $b=R_n$ gives
   \begin{align*}
     j(z_n) - j(u^\ast) &\leq   \frac{1}{2 \mu} \|g(z_n)\|^2 
        \leq \frac{1}{2\mu}\left(\frac{\sqrt{\rho}\mu}{L}\|g(u_n)\|^2+\frac{R_n^2}{1- \sqrt{\rho}}\right).
   \end{align*} %
   We again use \cref{eq:S1} to obtain the following bound for the full-step:   \begin{equation}\label{eq:full_step_bound}
        j(z_n) - j(u^\ast)\leq \sqrt{\rho} (j(u_n)-j(u^\ast))+\frac{R_n^2}{2 \mu(1- \sqrt{\rho})}.
   \end{equation}

\medskip\noindent\emph{Part 2: damped step $w_n$.} 
On $A_n^c$ with $\tau = \frac{m}{L}$, we denote $\bar{w}_n:= u_n - \tau H^{-1}_{S_n}(u_n)[g(u_n)]$
the damped step with exact gradient, such that $w_n = \bar{w}_n - \tau H^{-1}_{S_n}(u_n)[\pi_n]$. Then, by 
\eqref{eq:H_inv_bounds} and the definition of $\tau$, we have 
\begin{equation}\label{eq:w_n_error}
    \|w_n - \bar{w}_n\| \leq \frac{\tau}{m}\|\pi_n\| = \frac{\|\pi_n\|}{L}.
\end{equation}
By the descent lemma \eqref{eq:def_descent_lemma}, we have 
\begin{align*}
     j(\bar{w}_n) &\leq j(u_n) + \langle g(u_n) , \bar{w}_n - u_n \rangle + \frac{L}{2}\|\bar{w}_n - u_n\|^2 \\
    &=  j(u_n)- \tau \langle g(u_n) ,H_{S_n}^{-1}[g(u_n)]\rangle + \frac{L \tau^2}{2}\|H_{S_n}^{-1}[g(u_n)]\|^2.
\end{align*}
Applying the bound \eqref{eq:descent_lemma_hessian_bound} with $v = g(u_n)$ and $\|g(u_n)\|^2 \geq 2\mu(j(u_n)-j(u^*))$ from \eqref{eq:S1}, we have %
\begin{equation}\label{eq:exact_damped_step_contraction}
   j(\bar{w}_n) - j(u^\ast) \leq j(u_n)-j(u^*) - \frac{\tau}{2M}\|g(u_n)\|^2 \leq %
   \rho (j(u_n)-j(u^*)).
\end{equation}
The descent lemma \eqref{eq:def_descent_lemma} at $w_n, \bar{w}_n$ yields
\vspace{-1em}
\begin{align*}
    j(w_n)  &\leq j(\bar{w}_n) + \langle g(\bar{w}_n), w_n - \bar{w}_n\rangle + \frac{L}{2}\|w_n - \bar{w}_n\|^2\\
    &\leq j(\bar{w}_n)  + \frac{1}{L}\|g(\bar{w}_n)\| \|\pi_n\| + \frac{1}{2L}\|\pi_n\|^2\\
    &\leq j(\bar{w}_n) + \frac{1}{L}\left(\frac{\zeta}{2}\|g(\bar{w}_n)\|^2 + \frac{1}{2\zeta} \|\pi_n\|^2 \right) + \frac{1}{2L}\|\pi_n\|^2\\
    &\leq j(\bar{w}_n) + \zeta(j(\bar{w}_n)-j(u^*)) +  \frac{1+ \zeta^{-1}}{2L}\|\pi_n\|^2,
\end{align*}
where the second step uses the Cauchy--Schwarz inequality and \eqref{eq:w_n_error}, the third step  uses Young's inequality for $\zeta=(1/\sqrt{\rho}-1)$, and the final step applies \eqref{eq:S1}.
Combining with \eqref{eq:exact_damped_step_contraction}, we obtain
\vspace{-1em}
\begin{equation}\label{eq:damped_step_bound}
    j(w_n) - j(u^\ast) \leq \sqrt{\rho} (j(u_n)-j(u^\ast))+ \frac{\|\pi_n\|^2}{2L(1- \sqrt{\rho})}.
\end{equation}

\medskip\noindent\emph{Part 3: combining both steps.}
We now plug the full-step bound \eqref{eq:full_step_bound} and the damped step bound \eqref{eq:damped_step_bound} into \eqref{eq:split}. Using $\mathbf{1}_{A_n} + \mathbf{1}_{A_n^c} = 1$ gives 
\vspace{-1em}
\begin{equation}\label{eq:combined_pointwise}
   j(u_{n+1})-j(u^\ast) \leq \sqrt{\rho}(j(u_{n})-j(u^\ast)) + \frac{1}{2(1-\sqrt{\rho})}\left( \frac{\mathbf{1}_{A_n} R_n^2}{\mu} + \frac{\mathbf{1}_{A_n^c}\|\pi_n\|^2}{L}\right).
\end{equation}
We now take $\mathbb{E}[\cdot \mid \mathcal{F}_n]$ and bound the two noise terms by $V_n$ as defined in \cref{rem:iterate_noise_bound}. Further, the remaining indicator functions can be dropped because all summands are positive. 
To continue, note that $R_n^2 \leq 3\left( \delta_n^2 + \theta^2\|\pi_n\|^2 + \|\tilde{\pi}_n\|^2 \right)$. 
Further, we have  $\mathbb{E}[\|\pi_n\|^2 \mid \mathcal{F}_n] \leq \frac{V_n}{|X_n|}$ by \eqref{eq:variance_trick}. To bound $\|\tilde{\pi}_n\|^2$, we observe that $z_n$ is $\mathcal{F}_n^+$-measurable as defined in \eqref{eq:filtrations}, while $\tilde{X}_n$ is drawn independently of $\mathcal{F}_n^+$. Hence, we have with \cref{rem:conditional-variance-minibatch}, $|X_n| = |\tilde{X}_n|$, the tower property \cite[Theorem 4.52]{lord2014IntroductionComputationalStochastic}, and \cref{lem:initial_zn_bound}, that 
\[\mathbb{E}[\|\tilde{\pi}_n\|^2 \mid \mathcal{F}_n ] = \mathbb{E} \left[\mathbb{E}[\|\tilde{\pi}_n\|^2 \mid \mathcal{F}_n^+ ] \mid \mathcal{F}_n \right]  \leq \frac{\nu_0^2 + \nu_1^2 \mathbb{E}[\|z_n- u^\ast\|^2 \mid \mathcal{F}_n ]}{|X_n|} \leq \frac{C_z V_n}{|X_n|}.\]
Thus, 
\[\mathbb{E}[R_n^2 \mid \mathcal{F}_n ] \leq 3 \delta_n^2 + \frac{3(\theta^2 + C_z) V_n}{|X_n|},\]
which turns the conditional expectation of \eqref{eq:combined_pointwise} into 
\[\mathbb{E}[\Delta_{n+1} \mid \mathcal{F}_{n} ] \leq \sqrt{\rho} \Delta_n + \frac{3 \delta_n^2}{2 \mu(1-\sqrt{\rho})}  + \underbrace{\frac{1}{2(1-\sqrt{\rho})} \left(\frac{3 (\theta^2 + C_z)}{\mu}  + \frac{1}{L}\right)}_{=K} \frac{V_n}{|X_n|}.\]
Plugging in $V_n$ as defined in \eqref{eq:def_Vn} yields the explicit forms of $\rho_n, e_n.$ Lastly, with $|X_n| \geq X_{\min}$ as in \eqref{eq:X_min} we have 
\[\rho_n = \sqrt{\rho} + \frac{2 \nu_1^2 K}{\mu |X_n|} \leq \sqrt{\rho} +\frac{1- \sqrt{\rho}}{2} = \frac{1 + \sqrt{\rho}}{2}<1\]
which completes the proof. 
\end{proof}

We can now state global R-linear convergence of \cref{alg:globalized_SSN} in expectation.

\begin{theorem}[Global R-linear convergence in expectation]\label{thm:global_linear_combined}
    Let \cref{ass:lipschitz_mu_convex,ass:gradient_variance,ass:hessian_bounds} hold and $\{u_n\}$ be generated by \cref{alg:globalized_SSN} with $|{X}_n|=|\tilde{X}_n|\geq X_{\min}$ for all $n$. Let $\alpha \in (\bar{\rho},1)$ and suppose there exists $c_e >0$ such that $e_n \leq c_e \alpha^n$ for all $n \geq 0$. Then, defining $C:=  j(u_0)- j(u^\ast) + \frac{c_e}{\alpha - \bar{\rho}}$ , for all $n \geq 0$,
    \begin{align}
        \mathbb{E}[j(u_n)- j(u^\ast)] &\leq C \alpha^n\to 0,   \label{eq:e_n_geometric} \\
        \mathbb{E}[\|u_n - u^\ast\|^2 ] &\leq \frac{2}{\mu}C \alpha^n \to 0. \label{eq:u_n_quadratic_global_bound}
    \end{align}
    The condition on $e_n$ holds, for instance,
if $|X_n|^{-1} = \mathcal{O}(\alpha^n)$
and $\delta_n^2 = \mathcal{O}(\alpha^n)$.
\end{theorem}

\begin{proof}
    Taking the total expectation of \eqref{eq:global_linear} using $\rho_n \leq \bar{\rho}$ and $e_n \leq c_e\alpha^n$,
where $\rho_n, e_n$ are deterministic, yields $\mathbb{E}\big[j(u_{n+1}) - j(u^\ast)\big]
  \le \bar{\rho}\,\mathbb{E}\big[j(u_n) - j(u^\ast)\big]
    + c_e\,\alpha^n .$
    Unrolling the recursion with $\alpha > \bar{\rho}$ gives 
    \vspace{-1em}
    \begin{align*}
    \EE{j(u_n) - j(u^\ast)} &\leq \bar{\rho}^n (j(u_0) - j(u^\ast)) +c_e\sum_{k=0}^{n-1}\bar{\rho}^{n-1-k} \alpha^k\\
    &\leq \bar{\rho}^n (j(u_0) - j(u^\ast)) +\frac{c_e}{ \alpha- \bar{\rho}}\alpha^n \leq C \alpha^n.
\end{align*} Strong convexity \eqref{eq:S3} gives \eqref{eq:u_n_quadratic_global_bound}.
\end{proof}

\begin{remark}
    We note that no assumptions on Hessian batch sizes are needed for global R-linear convergence; i.e., $|S_n| \equiv 1$ suffices for global R-linear convergence. Hessian batch growth is needed only for the local, superlinear phase.
\end{remark}

\begin{remark}
    \label{rem:global_iterate_bounds}
    By \cref{thm:global_linear_combined}, we have $
        \sup_n \mathbb{E}[j(u_n)-j(u^\ast)] \leq C < \infty$
    and the noise level $V_n$ as defined in \cref{eq:def_Vn} is uniformly bounded in expectation, i.e.,
    \begin{equation}\label{eq:Vn_bound}
        \sup_n \mathbb{E}[V_n] \leq \bar{V} := \nu_0^2 + \frac{2 \nu_1^2}{\mu}C.
    \end{equation}
    With \cref{eq:variance_trick,eq:initial_zn_bound}, this implies
\begin{equation*}\label{eq:total_grad_var}
    \mathbb{E}\big[\|G_{X_n}(u_n) - g(u_n)\|^2\big] \leq \frac{\bar{V}}{|X_n|},
    \qquad
    \mathbb{E}\big[\|G_{\tilde{X}_n}(z_n) - g(z_n)\|^2\big] \leq \frac{C_z\bar{V}}{|{X}_n|}.
\end{equation*} 
\end{remark}

We conclude this subsection by showing almost sure convergence of the iterates and the objective gap.
\begin{lemma}[Almost sure convergence]\label{lem:almost_sure_convergence}
    Let the assumptions of \cref{thm:global_linear_combined} hold. Then, for the sequence of iterates $\{u_n\}$ generated by \cref{alg:globalized_SSN} it holds that 
    \begin{enumerate}[(i)]
        \item $\{\|u_n - u^\ast\|^2\}$ converges almost surely and $\lim_{n \to \infty} u_n = u^\ast$,
        \item $\{j(u_n) - j(u^\ast)\}$ converges almost surely and $\lim_{n \to \infty} j(u_n) = j(u^\ast)$.
    \end{enumerate}
\end{lemma}
\begin{proof}
    For every $\varepsilon>0$, Markov's inequality and \eqref{eq:u_n_quadratic_global_bound} give
\[
\sum_{n=1}^{\infty}
\mathbb{P}\bigl(\|u_n-u^*\|^2>\varepsilon\bigr) \leq \sum_{n=1}^{\infty} \frac{\EE{\|u_n-u^*\|^2}}{\varepsilon}
\leq
\frac{2C}{\mu\varepsilon}
\sum_{n=1}^{\infty}\alpha^n
<\infty.
\]
The first Borel--Cantelli lemma therefore yields
$\mathbb{P}\bigl(\|u_n-u^*\|^2>\varepsilon
\ \text{infinitely often}\bigr)=0.$
Since this holds for every $\varepsilon>0$, we conclude that $u_n$ converges to $u^*$ almost surely. In the same way, from \eqref{eq:e_n_geometric}, we have that $j(u_n)$ converges to $j(u^*)$ almost surely.
\end{proof}

\subsection{Switching to superlinear convergence}
We now show that the iterates will eventually reach a ball $B(u^\ast, r_\star)$ where the convergence behavior changes. Once inside the ball, all future full Newton steps pass the test \eqref{eq:ast}, remain inside the ball and converge R-superlinearly towards $u^\ast. $

\begin{assumption}[Local parameters]\label{ass:local_parameters}
    Let $\lambda:= \frac{\theta\mu}{L} \in (0,1)$
and fix an entry radius $r_\star > 0$ with $\omega(r_\star) < \lambda$, deterministic sequences $c_n \downarrow 0, \ \beta_n \downarrow 0$.
We assume the test slack satisfies 
\begin{equation}\label{eq:slack_lower_bound}
    \delta_n \geq \Big(\frac{L}{m} + 1 + \theta\Big)\, c_n
    \quad \text{for all sufficiently large } n.
  \end{equation}
We denote by $N_C$ the smallest index such that \eqref{eq:slack_lower_bound} holds and, for all $n \geq N_C$
\begin{equation}\label{eq:N_C_schedule}
  \omega(r_\star) + \frac{\beta_n}{m} \leq \lambda, \quad
  c_n \leq m(1-\lambda)\,r_\star.
\end{equation}
\end{assumption}
Since $\theta \in (0,1)$ and $\mu \leq L$, we have $\lambda \in (0,1)$  and $N_C$ exists since $\omega(r_\star) < \lambda$
and $c_n, \beta_n \downarrow 0$. Combined with the assumptions from \cref{thm:global_linear_combined}, the slack is thus bounded from above by $\delta_n^2 = \mathcal{O}(\alpha^n)$ and bounded from below by \eqref{eq:slack_lower_bound}.

The following lemma is an event-wise version of \cite[Lemma 2.3, Lemma 2.4]{bollapragada2019ExactInexactSubsampled}, with their Lipschitz-Hessian term replaced by the semismoothness modulus.

\begin{lemma}\label{lem:Eng_Enh}
Let \cref{ass:lipschitz_mu_convex,ass:newton_differentiability,ass:hessian_bounds} hold, fix $n \in \mathbb{N}$ and define, with $c_n, \beta_n$ from \cref{ass:local_parameters}, the gradient and Hessian noise events 
\begin{align*}
E_n^g&:= \left\{\|G_{X_n}(u_n) - g(u_n)\| \leq c_n \right\},\\
E_n^h&:= \left\{ \| (H_{S_n}(u_n) - h(u_n))[u_n - u^\ast]\| \leq \beta_n \|u_n - u^\ast\|\right\}.    
\end{align*}
Then, on $E_n^g\cap E_n^h$, the full step $z_n=u_n-H_{S_n}^{-1}(u_n)[G_{X_n}(u_n)]$ satisfies
\begin{equation}\label{eq:full_step_event}
    \|z_n - u^\ast\| \leq \left(\omega(\|u_n - u^\ast\|) + \frac{\beta_n}{m} \right)\|u_n - u^\ast\| + \frac{c_n}{m}.
\end{equation}

\end{lemma}
\begin{proof}
    Using the Hessian lower bound \cref{ass:hessian_bounds}, plugging $s = u_n-u^\ast$ into \eqref{eq:semismooth_modulus}, and then applying the definitions of $E_n^g, E_n^h$, we obtain
    \begin{align*}
        \|z_n - u^\ast\| &\leq \omega(\|u_n- u^\ast\|)\|u_n - u^\ast\| \\
        &\quad + \frac{1}{m}\|G_{X_n}(u_n) - g(u_n)\| + \frac{1}{m}\|(H_{S_n}(u_n) - h(u_n))[u_n - u^\ast]\|\\
        &\leq \left(\omega(\|u_n - u^\ast\|) + \frac{\beta_n}{m} \right)\|u_n - u^\ast\| + \frac{c_n}{m}.
    \end{align*}

\end{proof}
  
\begin{lemma}\label{lem:E_nt}
    Let \cref{ass:lipschitz_mu_convex,ass:newton_differentiability,ass:hessian_bounds,ass:local_parameters} hold.  In addition to $E_n^g,E_n^h$, we define the fresh test-batch event
\[
E_n^t:=\left\{ \|G_{\tilde X_n}(z_n)-g(z_n)\|\le c_n\right\},
\]
and set $E_n:=E_n^g\cap E_n^h\cap E_n^t$. If $n \geq N_C$ and $\|u_n - u^\ast\| \leq r_\star$, then the stochastic acceptance test \eqref{eq:ast} holds on the event $E_n$.
\end{lemma}

\begin{proof}
    We assume $E_n$ holds. On $E_n^t$, using \eqref{eq:S2}, we have 
    \[\|G_{\tilde{X}_n}(z_n)\| \leq  \|g(z_n)\| + \|G_{\tilde{X}_n}(z_n) - g(z_n)\|\leq \|g(z_n)\|+c_n\leq L\|z_n-u^\ast\|+c_n.\]
    By \cref{lem:Eng_Enh}, monotonicity of $\omega$, $\|u_n - u^\ast\| \leq r_\star$, \eqref{eq:N_C_schedule}, and $L\lambda = \theta \mu$, we have
    \[
  L\|z_n - u^\ast\|
  \leq L\left(\omega(r_\star) + \frac{\beta_n}{m}\right)
      \|u_n - u^\ast\|
    + \frac{L}{m}c_n
  \leq \theta\mu \|u_n - u^\ast\| + \frac{L}{m} c_n .
\]
On $E_n^g$, \eqref{eq:S2} yields $\mu\|u_n - u^\ast\| \leq \|g(u_n)\| \leq \|G_{X_n}(u_n) \| + c_n$. Combining the three estimates and $n \geq N_C$, we obtain 
\[\|G_{\tilde{X}_n}(z_n)\|  \leq \theta \|G_{X_n}(u_n)\| + \left(\frac{L}{m} + 1 + \theta \right)c_n \leq \theta\|G_{{X}_n}(u_n)\| + \delta_n \] 
which is the test \eqref{eq:ast}.
\end{proof}

\begin{lemma}[Bounded superlinearly shrinking radii]\label{lem:bounded_shrinking_Lambda}
Let \cref{ass:local_parameters} hold and let $N \geq N_C$. We define $\{\Lambda_n\}_{n \ge N}$ by
\begin{equation}\label{eq:def_Lambda}
    \Lambda_N:=r_\star,\qquad
\Lambda_{n+1}:=\left(\omega(\Lambda_n)+\frac{\beta_n}{m}\right)\Lambda_n+\frac{c_n}{m}.
\end{equation}
Then, the following hold:
\begin{enumerate}[(i)]
    \item $\Lambda_n \leq r_\star$ and 
    $\Lambda_{n+1}\leq \lambda \Lambda_n+\frac{c_n}{m}\le r_\star$ for all $n \geq N$,
    \item we have $\Lambda_n \to 0$ and hence $\omega(\Lambda_n) \to 0,$ and
    \item if $c_n \to 0$ Q-superlinearly, i.e.,~$\frac{c_n}{c_{n-1}}\to 0,$ then $\Lambda_n\to0$ Q-superlinearly, i.e.,~we have $\frac{\Lambda_{n+1}}{\Lambda_n} \to 0,$ as $n \to \infty.$
\end{enumerate}
    
\end{lemma}

\begin{proof}
    We show (i) by induction starting from $\Lambda_N = r_\star. $ If $\Lambda_n \leq r_\star$ for some $n \ge N \ge N_C$, the monotonicity of $\omega$ and \eqref{eq:N_C_schedule} give
     $\omega(\Lambda_n) + \frac{\beta_n}{m} \leq \omega(r_\star) + \frac{\beta_n}{m} \leq \lambda$. Hence,
    \[
        \Lambda_{n+1} \leq \lambda \Lambda_n + \frac{c_n}{m}\leq  \lambda r_\star + (1-\lambda) r_\star = r_\star.
    \]

    For (ii), we note that by (i), $0 \leq \Lambda_n \leq r_\star$ for all $n \geq N,$ which gives a finite $\bar{\Lambda}:= \limsup_{n \to \infty} \Lambda_n < \infty.$ Taking $\limsup$ on both sides of $\Lambda_{n+1} \leq \lambda \Lambda_n + \frac{c_n}{m}$ and using $c_n \to 0$ we have $\bar{\Lambda} \leq \lambda \bar{\Lambda}.$ 
    Since $\lambda < 1$ and $\Lambda_n \geq 0$ for all $n \in \mathbb{N}$, we get $\Lambda_n \to 0$, and $\omega(\Lambda_n) \to 0$ follows by \cref{ass:newton_differentiability} and \cref{rem:basic_SS}.
    
    To show (iii), we note that $\Lambda_N = r_\star > 0$ and, since all summands in \eqref{eq:def_Lambda} are nonnegative, $\Lambda_n \geq \frac{c_{n-1}}{m} > 0$ for $n > N$. Then, we divide \eqref{eq:def_Lambda} by $\Lambda_n>0$ to obtain 
    \begin{equation*}
        \frac{\Lambda_{n+1}}{\Lambda_n}=\omega(\Lambda_n)+\frac{\beta_n}{m}+\frac{c_n}{m\Lambda_n} \leq \omega(\Lambda_n) + \frac{\beta_n}{m} + \frac{c_n}{c_{n-1}} .
    \end{equation*}
  The first term tends to zero by (ii), the second by $\beta_n \to 0$ and the third by assumption. 
\end{proof}

We are now at a point where we can formulate the switch to superlinear convergence, i.e.,~that the iterates $\{u_n\}$ almost surely enter $B(u^\ast, r_\star)$, after which all future steps will be full stochastic semismooth Newton steps with $\tau_n = 1$ and the iterates converge R-superlinearly towards $u^\ast.$

The following theorem is reminiscent of \cite[Theorem 4.5]{milzarek2022local} which states that their method, with high probability, will switch to purely stochastic Newton steps as opposed to proximal gradient steps and that the iterates converge R-superlinearly for sufficiently quickly growing batch sizes. We show that this transition happens almost surely, after which almost surely only full Newton steps are accepted and we converge R-superlinearly. %

\begin{theorem}[Almost sure local R-superlinear convergence]\label{thm:local_r_superlinear}
    In addition to the assumptions of \cref{thm:global_linear_combined}, let  \cref{ass:newton_differentiability,ass:hessian_variance} and \cref{ass:local_parameters} hold. Fix $\{p_n\}_{n \in \mathbb{N}}$ with $p_n >0$ for all $n\in \mathbb{N}$ and $\sum_{n \ge 0} p_n < \infty$,
    and impose the additional batch floors
    \begin{equation}\label{eq:as_batch_floors}
      |X_n| \geq \frac{C_z\,\bar{V}}{c_n^2\, p_n},
      \qquad
      |S_n| \geq \frac{\sigma^2}{\beta_n^2\, p_n}
      \qquad \text{for all } n \ge 0 .
    \end{equation}
    Then there exists an almost surely finite random index $N^\ast \geq N_C$ such that almost surely, for every  $n \ge N^\ast$:
    \begin{enumerate}[(i)]
        \item the test \eqref{eq:ast} is accepted, so $\tau_n = 1$ and
        $u_{n+1} = z_n$; in particular, at most $N^\ast$
        damped fallback steps are taken,
        \item $\|u_n - u^\ast\| \leq  r_\star$, so
        the iterates never leave $B(u^\ast, r_\star)$,
        \item if $c_n \to 0$ at a Q-superlinear rate, then $u_n \to u^\ast$ R-superlinearly almost surely.
    \end{enumerate}
    
\end{theorem}

\begin{proof}
\textit{Step 1 (almost surely, only finitely many $E_n^c$ occur).}
We recall the noise events $E_n := E_n^g \cap E_n^h \cap E_n^t$ from \cref{lem:Eng_Enh} and \cref{lem:E_nt} with the convention that $E_n^h$ holds on $\{u_n = u^\ast\}$. Applying Markov's inequality to the conditional $\mathbb{P}((E_n^g)^c)$, with both sides of $(E_n^g)^c$ squared and using \eqref{eq:variance_trick}, we obtain 
\[\mathbb{ P}\left((E_n^g)^c \mid \mathcal{F}_n\right)\leq\frac{\mathbb{E}[\|G_{X_n}(u_n)-g(u_n)\|^2 | \mathcal{F}_n]}{c_n^2}\leq \frac{\nu_0^2+\nu_1^2\|u_n- u^\ast\|^2}{c_n^2|X_n|} \leq \frac{V_n}{c_n^2|X_n|}.\]
Taking total expectation and using \eqref{eq:Vn_bound} the floor
\eqref{eq:as_batch_floors} and $C_z \geq 1$, we obtain
\[\mathbb{ P}\left((E_n^g)^c\right) = \mathbb{E}\left[ \mathbb{ P}\left((E_n^g)^c \mid \mathcal{F}_n\right)\right]\leq \frac{\bar{V}}{c_n^2|X_n|} \leq \frac{\bar{V}C_z}{c_n^2|X_n|} \leq p_n.\]
The same follows analogously for $(E_n^t)^c$, where we condition on $\mathcal{F}_n^+$ as defined in \eqref{eq:filtrations} which already has information of $u_n, z_n$. However, $\mathcal{F}_n^+$ is independent from $\tilde{X}_n$, which gives 
\[\mathbb{ P}\left((E_n^t)^c\right) = \mathbb{E}\left[ \mathbb{ P}\left((E_n^t)^c \mid \mathcal{F}_n^+\right)\right] \leq p_n.\]

For the Hessian noise event, on the $\mathcal{F}_n$-measurable event
$\{u_n \neq u^\ast\}$, we can proceed analogously with \cref{ass:hessian_variance} to obtain 
\[\mathbb{P}\left((E_n^h)^c \mid \mathcal{F}_n\right)\leq\frac{\mathbb{E}[\|(H_{S_n}(u_n)-h(u_n))[u_n - u^\ast]\|^2 \mid \mathcal{F}_n]}{\beta_n^2\|u_n - u^\ast\|^2}\leq\frac{\sigma^2}{\beta_n^2|S_n|}\leq p_n,\] and hence 
$\mathbb{P}\left((E_n^h)^c\right) \leq p_n,$ since on $\{u_n = u^\ast\}$ the event $E_n^h$ holds by convention. 

By subadditivity,
\[\mathbb{P}(E_n^c)\leq\mathbb{P}((E_n^g)^c)+\mathbb{P}((E_n^h)^c)+\mathbb{P}((E_n^t)^c)\leq p_n+p_n+p_n=3 p_n \]
for every $n\geq 0$ and therefore 
\[\sum_{n \geq 0} \mathbb{P}(E_n^c) \leq 3 \sum_{n\geq 0} p_n  < \infty.\]
By the Borel-Cantelli lemma, $\mathbb{P}( \limsup_{n\to \infty} E_n^c ) = 0$, i.e., the probability of $E_n^c$ happening infinitely often is zero. Thus, there exists an almost surely finite random index 
\[N_E := \inf\{N \in \mathbb{N} : \mathbf{1}_{E_n} = 1 \ \text{ for all } n \geq N\}.\]

\textit{Step 2 (almost surely, the iterates eventually remain in $B(u^\ast, r_\star)$).}
By \cref{lem:almost_sure_convergence}, we have $u_n \to u^\ast$ almost surely, so almost surely $\|u_n - u^\ast\| > r_\star$ for at most finitely many $n$. Hence
\[N_B := \inf\{N \in \mathbb{N} : \|u_n - u^\ast\| \leq r_\star \ \text{ for all } n \geq N\}\]
is almost surely finite.

\textit{Step 3 (almost surely, the iterates transition to R-superlinear convergence).} We set 
$N^\ast := \max(N_E, N_B, N_C)$, which is almost surely finite and let $\{\Lambda_n\}_{n \geq N^\ast}$ be the sequence from \cref{lem:bounded_shrinking_Lambda} with $N = N^\ast$.
We now show by induction that for all $n \geq N^\ast,$ the test \eqref{eq:ast} is accepted and we have $\|u_n - u^\ast\| \leq \Lambda_n \leq r_\star$. For the base case $n = N^\ast$  we have  $\|u_{N^\ast} - u^\ast\| \leq r_\star = \Lambda_{N^\ast}$ since $N^\ast\geq N_B.$ We now assume $\|u_n - u^\ast\| \leq \Lambda_n \leq r_\star$ for some $n \geq N^\ast.$ 
First, since $n \geq N_E$, the event $E_n$ holds, and since $n \geq N_C$ and $\|u_n - u^\ast\| \leq r_\star$,
\cref{lem:E_nt} shows that the test \eqref{eq:ast} is accepted,
i.e., $u_{n+1} = z_n$.
Further, \cref{lem:Eng_Enh} and the definition of $\Lambda_n$ \eqref{eq:def_Lambda} along with the monotonicity of $\omega$ gives
\begin{align*}
    \|u_{n+1} - u^\ast\| &\leq \left(\omega(\|u_n - u^\ast\|) + \frac{\beta_n}{m} \right)\|u_n - u^\ast\| + \frac{c_n}{m}\\
    & \leq \left(\omega(\Lambda_n) + \frac{\beta_n}{m} \right)\Lambda_n + \frac{c_n}{m} = \Lambda_{n+1} \leq r_\star.
\end{align*}

This finishes the induction and proves (i) and (ii).  Lastly, (iii) is shown by \cref{lem:bounded_shrinking_Lambda}.  
\end{proof}

\begin{remark}
    We briefly state what choices for the batch sizes $|X_n|, |S_n|$ suffice for the global and local phase. 

    \paragraph{Global phase} Geometric growth in the gradient batches
    $|X_n| = \lceil \max(1, X_{\min})  q^n \rceil$ with $q > 1$ and slack $\delta_n = \delta_0|X_n|^{-1/2}$ for $\delta_0 > 0$ as well as a constant Hessian batch size $|S_n| = 1$ yield $e_n \leq c_e \alpha^n $ with $\alpha := \max(q^{-1}, \alpha_0)$ for $\alpha_0 \in (\bar{\rho}, 1).$ This 
    satisfies the assumptions of \cref{thm:global_linear_combined}. 

    \paragraph{Local phase} To fix the floors \eqref{eq:as_batch_floors} we first set $p_n = (n+1)^{-2}$ such that $\sum_{n}p_n < \infty$. Then, we can set
    $c_n^2 = (n+1)^2\frac{C_z \bar{V}}{|X_n|}.$ Thus, to obtain $c_n \to 0$ at a Q-superlinear rate, we need $|X_n|^{-1}\to 0$ at a Q-superlinear rate which suggests batch schedules such as 
    \[|X_n| = \lceil X_{\min}\,q^{n(n+1)/2}\rceil  \quad \text{or} \quad |X_n| = \lceil X_{\min}\,q^{n \ln(n+1)}\rceil \]
    for some $q >1$. For the generalized Hessian batches, we fix $\beta_n = \beta_0/(n+1)$ for some $\beta_0 >0$ such that $\beta_n \downarrow 0$, and we obtain 
    \[|S_n| = \lceil\sigma^2(n+1)^4/\beta_0^2\rceil.\]
    Lastly, as fixed in \cref{ass:local_parameters}, we set 
 $\delta_n = \left(\frac{L}{m} + 1 + \theta \right)c_n.$
    
\end{remark}

\begin{remark}[Oracle complexity: PDE solves versus iterations]\label{rem:oracle_complexity}
    One gradient sample $G_{X_n}(u_n)$ requires one state and one adjoint solve, and each application of $H_{S_n}(u_n)$ inside an iterative solver such as CG requires one linearized state and  adjoint solve per sample. Single sample SGD $u_{n+1} = u_n - \tau_n G_{\xi^n}(u_n)$ with optimal step sizes $\tau_n>0$
    reaches $\mathbb{E}[\|u_n - u^\ast\|^2] \leq \varepsilon$ after $\mathcal{O}(\varepsilon^{-1})$ iterations (see \cite{bottou2018optimization,nemirovsky1983problem}).
Now consider \cref{alg:globalized_SSN} under \cref{thm:global_linear_combined} with 
$|X_n| = |\tilde X_n| = \lceil X_{\min}\alpha^{-n}\rceil$, $\alpha \in (\bar\rho, 1)$, $\delta_n^2 = \mathcal{O}(\alpha^n)$ and $|S_n| \equiv 1$. Then, 
\eqref{eq:u_n_quadratic_global_bound} gives $\varepsilon$-accuracy after 
$n_\varepsilon = \left\lceil \frac{\log(2C/  \mu\varepsilon)}{\log(1/\alpha)}\right\rceil = \mathcal{O}(\log(\varepsilon^{-1}))$
 iterations, at the cost of 
$4\sum_{n=0}^{n_\varepsilon}|X_n|
\leq \frac{8CX_{\min}}{\alpha(1-\alpha)\mu\varepsilon}
+ 4(n_\varepsilon + 1)
= \mathcal{O}(\varepsilon^{-1})$ 
 PDE solves for the gradients and $\mathcal{O}(\log(\varepsilon^{-1}))$ Hessian applications. Summarizing, this leaves us at $\mathcal{O}(\varepsilon^{-1})$ state solves for the global SSN phase, the same as SGD, but compressed into $\mathcal{O}(\log(\varepsilon^{-1}))$ iterations. 
We note that for the local superlinear phase it is hard to derive meaningful complexity results since we do not have an explicit convergence rate, as the semismoothness modulus $\omega$ is not known explicitly. 
\end{remark}

\section{Numerical examples}\label{sec:numerics}
We apply \cref{alg:globalized_SSN} to PDE-constrained optimization problems and first fix relevant notation.
Let $D \subset \mathbb{R}^d$, $d \leq 3$ be a bounded Lipschitz domain. For $p \in [1,\infty]$, the space $L^p(D)$ denotes the usual Lebesgue space with the Lebesgue measure on $D$ \cite[Chapter 2]{adams2003SobolevSpaces}.
Further, we denote by $W^{s,p}(D)$ the Sobolev space of functions in $L^p(D)$ with weak derivatives up to order $s$ belonging to $L^p(D)$ \cite[Chapter 3]{adams2003SobolevSpaces}. For $p=2$, we write $H^s(D) := W^{s,2}(D)$, and by $H^1_0(D)$ we denote the closure of $C_c^\infty(D)$ in $H^1(D)$ \cite[Chapter 3]{adams2003SobolevSpaces}. 
For Banach spaces $X,Y,$ we write $\mathcal{L}(X,Y)$ for the space of bounded linear operators from $X$ to $Y$ and $X\hookrightarrow Y$ for a continuous embedding. In the following, for a measurable set $E \subseteq D,$ the notation $\operatorname{meas}(E)$ refers to its Lebesgue measure and we write $\mathbf{1}_E$ for its characteristic function which is to be understood pointwise.
Finally, for $v: \mathbb{R}^d \to \mathbb{R}$ and $w: \mathbb{R}^d \to \mathbb{R}^d$ we denote by $\nabla v(x)= ( \partial v(x)/\partial x_1, \dots, \partial v(x)/\partial x_d)^\top$ and $\nabla \cdot w(x) = \partial w_1(x)/ \partial x_1 + \dots + \partial w_d(x) /\partial x_d$ the gradient and divergence, respectively. Both only act on the spatial variable $x$.

\subsection{Linear elliptic PDE with Huber loss} \label{sec:linear_elliptic}
Compared to the classical $L^2(D)$-tracking, an $L^1(D)$-tracking term is more robust to outliers \cite{guitton2003robust, huber1992robust, tavolato2015use}, but it is not differentiable. Smoothing it on a small neighborhood around the origin yields the Huber loss, which acts as $L^2(D)$-tracking on the smoothing neighborhood and as $L^1(D)$-tracking otherwise. Hence, it is an element of $SC^1$ and, thus, suitable for our purposes. Further, it is also useful to demonstrate the need for a globalization strategy.

\subsubsection{Setup}
 We set $D:= (0,1)^2$. We wish to minimize the cost functional
 \begin{equation}\label{eq:numerics_elliptic_objective}
     \min_{u \in L^2(D)} j(u) :=  \EE{J(u, \xi)} :=\EE{ \int_D \Psi_r(y(x,\xi)-y_d(x)) \, \mathrm{d}x} + \frac{\mu}{2}\|u\|_{L^2(D)}^2
 \end{equation}
 over the Hilbert space $U:= L^2(D)$ subject to the random linear elliptic PDE
\begin{equation}\label{eq:example_state}
    \begin{aligned}
        -\nabla \cdot \bigl( a(x, \xi)\, \nabla y(x, \xi) \bigr) &= u(x) + f(x, \xi), && x \in D,\ \xi \in \Xi, \\
y(x, \xi) &= 0, && x \in \partial D.
    \end{aligned}
\end{equation}
Above, we have $\mu> 0$ and the Huber loss function and its derivative
\[\Psi_r (t) := \begin{cases}
    \frac{1}{2}t^2, &|t|\leq r,\\
    r |t| - \frac{1}{2}r^2, &|t|>r,
\end{cases}\qquad \Psi_r'(t) = \min(\max(t, -r), r), \] for $r >0.$
Both are understood as superposition operators, i.e., for measurable $v\colon D \to \mathbb{R}$, we set $\Psi_r(v)(x)= \Psi_r(v(x))$ for a.e.~$x\in D$ and analogously for $\Psi_r'.$ 
For the noisy diffusion parameter $a: D \times \Xi \to \mathbb{R}$, we use the truncated Karhunen--Loève expansion 
\[a(x, \xi)= 1 + \sum_{k=1}^4 \vartheta_k \xi^l_k \psi_k(x),\] with i.i.d.~$\xi^l_k\, {\sim}\, \mathcal{U}(-1,1)$, i.e., uniformly distributed on $(-1,1)$ \cite[Example 4.9]{lord2014IntroductionComputationalStochastic},
$\vartheta = (0.25,\ 0.15,\ 0.10,\ 0.05)$ and \[
\psi_1 = \cos(\pi x_1),\ \psi_2 = \cos(\pi x_2),\ \psi_3 = \cos(\pi x_1)\cos(\pi x_2),\ \psi_4 = \cos(2\pi x_1)\cos(\pi x_2).\] This gives 
$a(x,\xi) \in (a_{\min}, a_{\max}) = (0.45,\ 1.55)$ for every $\xi,$ so uniform ellipticity holds. For the right-hand side, we set \[f(x,\xi) = 0.5\,\big(\xi^r_1\, \sin(\pi x_1)\sin(\pi x_2) + \xi^r_2\, \sin(2\pi x_1)\sin(\pi x_2)\big)\]
 with i.i.d.~$\xi^r_k  \, \sim \,  \mathcal{N}(0,1),$ i.e., standard normally distributed \cite[Example 4.10]{lord2014IntroductionComputationalStochastic}. Thus, $ \xi = (\xi_1^l, \xi_2^l, \xi_3^l, \xi_4^l, \xi_1^r, \xi_2^r) \in \Xi \subset \mathbb{R}^6$.
For the optimization parameters, we set 
\[y_d(x) = 0.15\, \sin(\pi x_1)\sin(\pi x_2), \quad r = 0.2, \quad \text{and }\mu = 10^{-4}.\]

 A necessary and sufficient condition for $u^*$ to be a solution to the problem \eqref{eq:numerics_elliptic_objective} is given by the Riesz gradient representation
\begin{equation*}
    g(u) = \nabla j(u) = \EE{p(\cdot, \xi)} + \mu u =0,
\end{equation*}
where for each $\xi \in \Xi$, $p(\cdot, \xi)\in H^1_0(D) \hookrightarrow L^2(D)$ is the adjoint state determined by 
\begin{equation}\label{eq:example_adjoint}
    \begin{aligned}
         -\nabla \cdot \bigl( a(x, \xi)\, \nabla p(x, \xi) \bigr) &= \Psi'_r\left(y(x, \xi) - y_d(x)\right)  , && x \in D \\
        p(x, \xi) &= 0, && x \in \partial D.
    \end{aligned}
\end{equation}
See, e.g.,~\cite[Section 2.8]{troltzsch2010optimal} for the deterministic setting.

\subsubsection{Verification of assumptions}
\paragraph{Assumption 1: strong convexity and Lipschitz gradient}
Immediately, $\mu$-strong convexity follows by the quadratic Tikhonov term $\frac{\mu}{2}\|u\|_{L^2(D)}^2$ with $\mu >0.$ 
Since $a(x, \xi) \in [a_{\min}, a_{\max}]$ for every $\xi \in \Xi$, the associated bilinear form \[\mathfrak{a}_\xi(y,v):=\int_D a(x, \xi) \nabla y(x)\cdot\nabla v(x) \, \mathrm{d}x, \quad y,v \in H^1_0(D)\] is bounded, symmetric and coercive independently of $\xi.$
By the Lax--Milgram lemma \cite[Corollary 5.8]{brezis2011FunctionalAnalysisSobolev}, for each $\xi \in \Xi$, the solution operator associated to the state equation \eqref{eq:example_state} $S(\xi)\colon L^2(D) \to H^1_0(D)$ is well-defined. Since $H^1_0(D) \hookrightarrow L^2(D)$, we consider $S(\xi)$ as an operator on $L^2(D)$ whenever convenient and suppress the embedding. Symmetry of $\mathfrak{a}_\xi$ makes $S(\xi)$ self-adjoint on $L^2(D)$, and by Poincaré's inequality \cite[Corollary 9.19]{brezis2011FunctionalAnalysisSobolev} there exists $0<\bar{S}<\infty$ such that we have a $\xi$-independent bound
\[\|S(\xi)\|_{L^2(D) \to L^2(D)} \leq \bar{S} \qquad \text{for all } \xi \in \Xi.\]
We can further write 
\[y(\cdot, \xi) = S(\xi)(u + f(\cdot, \xi)).\]
Since $S(\xi)$ is self-adjoint, we can express the adjoint equation similarly as 
\[p(\cdot, \xi) = S(\xi)\Psi'_r(y(\cdot, \xi) -y_d).\]
Thus, for two controls $u_1, u_2$ and respective states and adjoints $y_1, y_2, p_1, p_2$, we obtain 
\begin{align*}
      \|p_1(\cdot,\xi) - p_2(\cdot,\xi)\|_{L^2(D)} &\leq \bar{S}\|\Psi_r'(y_1(\cdot, \xi)-y_d) - \Psi'_r(y_2(\cdot, \xi) - y_d)\|_{L^2(D)}\\
      &\leq \bar{S}\|y_1(\cdot, \xi)-y_2(\cdot, \xi)\|_{L^2(D)} \leq \bar{S}^2\|u_1 - u_2\|_{L^2(D)},
\end{align*}
where we use that $\Psi_r'=\min(\max(\cdot,-r),r)$ is the (pointwise) projection onto $[-r,r]$ \cite[Theorem 2.28]{troltzsch2010optimal} and is therefore nonexpansive \cite[Proposition 5.3]{brezis2011FunctionalAnalysisSobolev}.
Thus, we have
\begin{align*}
    \|g(u_1)- g(u_2)\|_{L^2(D)} &\leq\| \EE{p_1(\cdot, \xi)}-  \EE{p_2(\cdot, \xi)} \|_{L^2(D)} + \mu\|u_1 - u_2\|_{L^2(D)}\\ &\leq (\bar{S}^2 + \mu) \|u_1 -u_2\|_{L^2(D)}.
\end{align*}

\paragraph{Assumption 2: Newton differentiability of the gradient}
Since we are working with the nonsmooth $\min$ and $\max$ operators, Newton differentiability needs to be treated with care. The operators $\min, \max \colon L^q(D) \to L^p(D)$ are only Newton differentiable for $1 \leq p < q <\infty$ (cf. \cite[Theorem 4.5]{hintermuller2010SemismoothNewtonMethods}, \cite[Proposition 1.1]{ito2003SemiSmoothNewton}). 
Thus, when considering $\min\colon L^q(D) \to L^2(D)$ (respectively $\max$) as it appears in \eqref{eq:example_adjoint}, we need $q>2.$ For any $\xi \in \Xi$, the control-to-state map corresponding to \eqref{eq:example_state} $u \mapsto y(\cdot, \xi)$ belongs to $\mathcal{L}(L^2(D), H^1_0(D))$. 
Further, $H^1_0(D) \hookrightarrow L^q(D)$ with $q\in [2, \infty]$ for $d=1,$ $q\in [2, \infty)$ for $d=2$, and $q\in [2, \frac{2d}{d-2}]$ for $d\geq 3$ \cite[Corollary 9.14]{brezis2011FunctionalAnalysisSobolev}. In any case, we find a $q>2$
and the arguments for $\min(\max(\cdot, -r), r)$ are compatible with its Newton differentiability. 
 Thus, for $z \in L^q(D)$, a generalized derivative of $\Psi_r'(z) = \min(\max(z, -r), r)$, viewed as a mapping from $L^q(D)$ to $L^2(D)$, is given by multiplication with
$\mathbf{1}_{\{|z| \leq r\}}$. Therefore, we have a generalized derivative $\delta_p(\cdot, \xi) \in L^2(D)$ of $p(\cdot, \xi)$ given by the linearized state and adjoint equations
\begin{equation*}
    \begin{aligned}
     -\nabla \cdot \bigl( a(x, \xi)\, \nabla \delta_y(x, \xi) \bigr) &= \delta_u(x), && x \in D,\ \xi \in \Xi, \\
        \delta_y(x, \xi) &= 0, && x \in \partial D,\\
         -\nabla \cdot \bigl( a(x, \xi)\, \nabla \delta_{p}(x, \xi) \bigr) &=\mathbf{1}_{\{|y(\cdot, \xi) - y_d| \leq r\}} \delta_y(x, \xi)  , && x \in D,\ \xi \in \Xi\\
         \delta_p(x, \xi) &= 0, && x \in \partial D,
   \end{aligned}
\end{equation*}
for $\delta_u \in L^2(D)$ with $\delta_y \in H^1_0(D) \hookrightarrow L^2(D).$
Therefore,
\begin{equation}\label{eq:example_generalized_hessian}
    h(u)[\delta_u] = \mathbb{E}[\delta_p(\cdot, \xi)] + \mu \delta_u = \EE{S(\xi)\mathbf{1}_{\{|y(\cdot, \xi) - y_d| \leq r\}}
      \,S(\xi)\delta_u} + \mu\delta_u
\end{equation}
is a generalized Hessian (see also \cite[Section 33]{clason2026introduction} for the deterministic computation) where integrability holds since $S(\xi)$ is bounded for all $\xi \in \Xi$.

\paragraph{Assumption 3,4: variance bounds for gradient and generalized Hessian}
For both bounds, we first note that for every $L^2(D)$-valued random variable $X$ with finite second moment, we have 
\begin{equation}\label{eq:variance_inequality}
    \EE{\|X- \EE{X}\|_{L^2(D)}^2} = \EE{\|X\|_{L^2(D)}^2} - \|\EE{X}\|_{L^2(D)}^2 \leq  \EE{\|X\|_{L^2(D)}^2} .
\end{equation}
Further, we note that for every $u \in U, \ \xi \in \Xi$, we have the unbiased gradient estimator
\[\mathbb{E}[G_{\xi}(u)] = \mathbb{E}[p(\cdot, \xi) + \mu u ] = \mathbb{E}[p(\cdot, \xi) ]+ \mu u  = g(u).\] Then, using \eqref{eq:variance_inequality} with $X= p(\cdot, \xi)$, we have
\begin{align*}
    \mathbb{E}[\|G_{\xi}(u) - g(u)\|_{L^2(D)}^2] 
    &= \mathbb{E}[\|p(\cdot, \xi)  - \mathbb{E}[p(\cdot, \xi)] \|_{L^2(D)}^2]\\
    &\leq \mathbb{E}[\|p(\cdot, \xi)  \|_{L^2(D)}^2]\\
    &\leq \bar{S}^2 \mathbb{E}[\| \Psi'_r(y(\cdot, \xi) - y_d) \|_{L^2(D)}^2]\\
    &\leq\bar{S}^2  r^2 \operatorname{meas}(D)=: \nu_0^2,
\end{align*}
where we used 
\[\|\Psi'_r(y(x, \xi) - y_d(x))\|_{L^2(D)}^2 = \int_D |\Psi_r'(y(x, \xi) - y_d(x))|^2 \, \mathrm{d} x 
     \leq  r^2 \int_D 1 \, \mathrm{d} x = r^2\operatorname{meas}(D).\]
In this example, we  even have $\nu_1 = 0$, so uniformly bounded gradient variance without affine scaling. 
Unbiasedness and bounded variance for the generalized Hessian with $\sigma^2 = \bar{S}^4$ is analogous. 

\paragraph{Assumption 5: Hessian bounds}
For $\zeta \in \Xi$, using self-adjointness of $S(\zeta)$ we have
\begin{align*}
    \langle H_{\zeta}(u)[v], v\rangle_{L^2(D)} &= \langle S(\zeta)\left( \mathbf{1}_{\{|y(\cdot, \zeta) - y_d| \leq r\}}S(\zeta) v\right), v \rangle_{L^2(D)} + \mu \|v\|_{L^2(D)}^2 \\
    &= \int_D \mathbf{1}_{\{|y(\cdot, \zeta) - y_d| \leq r\}} |S(\zeta) v |^2 \, \mathrm{d} x + \mu\|v\|_{L^2(D)}^2,
\end{align*}
which gives the bounds $m = \mu$ and $M= \bar{S}^2 + \mu = L$.

\subsubsection{Numerical results}
For the numerical experiments, we use a uniform triangulation of piecewise linear, globally continuous finite elements with 121 nodes with mesh size $h_{\min}\approx 0.14$ which was built using FEniCS \cite{langtangen2016SolvingPDEsPython}. The Newton systems $H_{S_n}(u_n)[d_n] = -G_{X_n}(u_n)$ are solved using a matrix free Krylov method, specifically the conjugate gradient (CG) method with relative tolerance $10^{-3}.$ For the problem specific parameters, we set $m = \mu.$ To approximate $L$ (respectively $M$), in the generalized Hessian $h(u)[\delta_u]$, we replace the generalized derivative of the nonsmooth term $\Psi'_r$ by its pointwise upper bound $1$ which leaves the resulting mapping
\[\delta_u \mapsto \EE{S(\xi)S(\xi)}\delta_u + \mu\delta_u.\] For the discretized counterpart, we can then estimate its largest eigenvalue by using e.g.~\texttt{scipy.linalg.eigh} in the $L^2(D)$ inner product for a large batch of $\xi \in \Xi$ and picking the maximum and set this to $M=L$.

To test the globalization strategy, we initialize with $u_0 \equiv 10$ such that the corresponding state is sufficiently far away from $y_d$ and we are (for  most interior points $x \in D$) in the purely $L^1(D)$-tracking regime.  We run \cref{alg:globalized_SSN} and track convergence in terms of the norm $\|u_n - u^\ast_{\text{ref}}\|_{L^2(D)}$, where $u^\ast_{\text{ref}}$ was computed a priori using sample average approximation (SAA) with a fixed batch $\{\xi_i\}_{i= 1}^{N}$ with $N= 256,000$. Since the expected Monte Carlo accuracy $\|u^\ast - u^\ast_{\text{ref}}\|_{L^2(D)} \approx 1/\sqrt{N} \approx 2\cdot 10^{-3}$, we cannot expect to track convergence beyond this error. 
We will compare \cref{alg:globalized_SSN} with full-step SSN, i.e., SSN with a fixed step size $\tau_n = 1$. Further, we compare the results to stochastic gradient descent (SGD) with fixed sample size $|X_n|=1$ and step size $\tau_n = \frac{1}{(n+1)\mu}$ as a baseline (see e.g.~\cite{bottou2018optimization}). 
To test \cref{alg:globalized_SSN}, we differentiate between a global and a local batch size strategy in line with the differing requirements of \cref{thm:global_linear_combined} and \cref{thm:local_r_superlinear} (see \cref{tab:schedules_elliptic}).
\begin{table}[htbp]
  \centering
  \caption{Batch sizes for the elliptic problem. The switch from global to local phase happens at iteration $n = n_s$.}
  \label{tab:schedules_elliptic}
  \begin{tabular}{llcc}
    \toprule
    Phase & Iterations
      & $|X_n| = |\tilde X_n|$ & $|S_n|$ \\
    \midrule
    global & $n < n_s$
      & $2^n$
      & $1$ \\[4pt]
    local & $n \geq n_s$, $k := n - n_s$
      & $\max\Big( 2^n,\,
          \Big\lceil \dfrac{c_0}{c_k^2\, p_k } \Big\rceil \Big)$
      & $\Big\lceil \dfrac{\beta_0}{\beta_k^2\, p_k} \Big\rceil$ \\
    \bottomrule
  \end{tabular}
\end{table}
The batch sizes in \cref{tab:schedules_elliptic} are built from
the sequences of \cref{thm:local_r_superlinear} with 
\begin{equation}\label{elliptic_sequences}
    q = 0.8, \qquad c_n = q^{n(n+1)}, \qquad \beta_n  = q^n,   \qquad
  p_n = (n+1)^{-1.1}.
\end{equation}
This satisfies $c_n / c_{n-1} = q^{2n} \to 0$ (Q-superlinear convergence), $\beta_n \downarrow 0$ and $\sum_n p_n < \infty.$
The scalings $c_0, \beta_0$ are chosen such that the local batch sizes
connect to the global ones at $n = n_s$.
The switch from the global to the local schedule at iteration $n = n_s$ is heuristic. Since the local batch sizes in \cref{tab:schedules_elliptic} never fall below the global ones, the switch does not affect the assumptions of \cref{thm:global_linear_combined}. The switch only decides when to start paying for the larger batches required by \cref{thm:local_r_superlinear}. Switching too early wastes samples on damped steps far from $u^\ast$ while switching too late delays the superlinear phase.
Once the test \eqref{eq:ast} has been accepted in two consecutive iterations, we switch from the global to the local phase. As a safeguard, we switch at the latest at $n = 5.$ While this is heuristic, it proved successful in numerical experiments. 

For the full-step SSN method, we stick with the global scheme, since the method will not converge if initialized too far away from the solution. We run both semismooth Newton methods for 9 iterations and track the amount of PDE solves used throughout, i.e., for all gradient evaluations, tests and CG Hessian-vector products. We set the maximum number of iterations for SGD to half the number of total PDE solves of \cref{alg:globalized_SSN}, since SGD requires two PDE solves per iteration. This leaves \cref{alg:globalized_SSN} and SGD with the same number of PDE solves for a fair comparison.

\begin{figure}[ht]
    \centering
    \begin{subfigure}{0.48\textwidth}
        \centering
            \includegraphics[width=\linewidth]{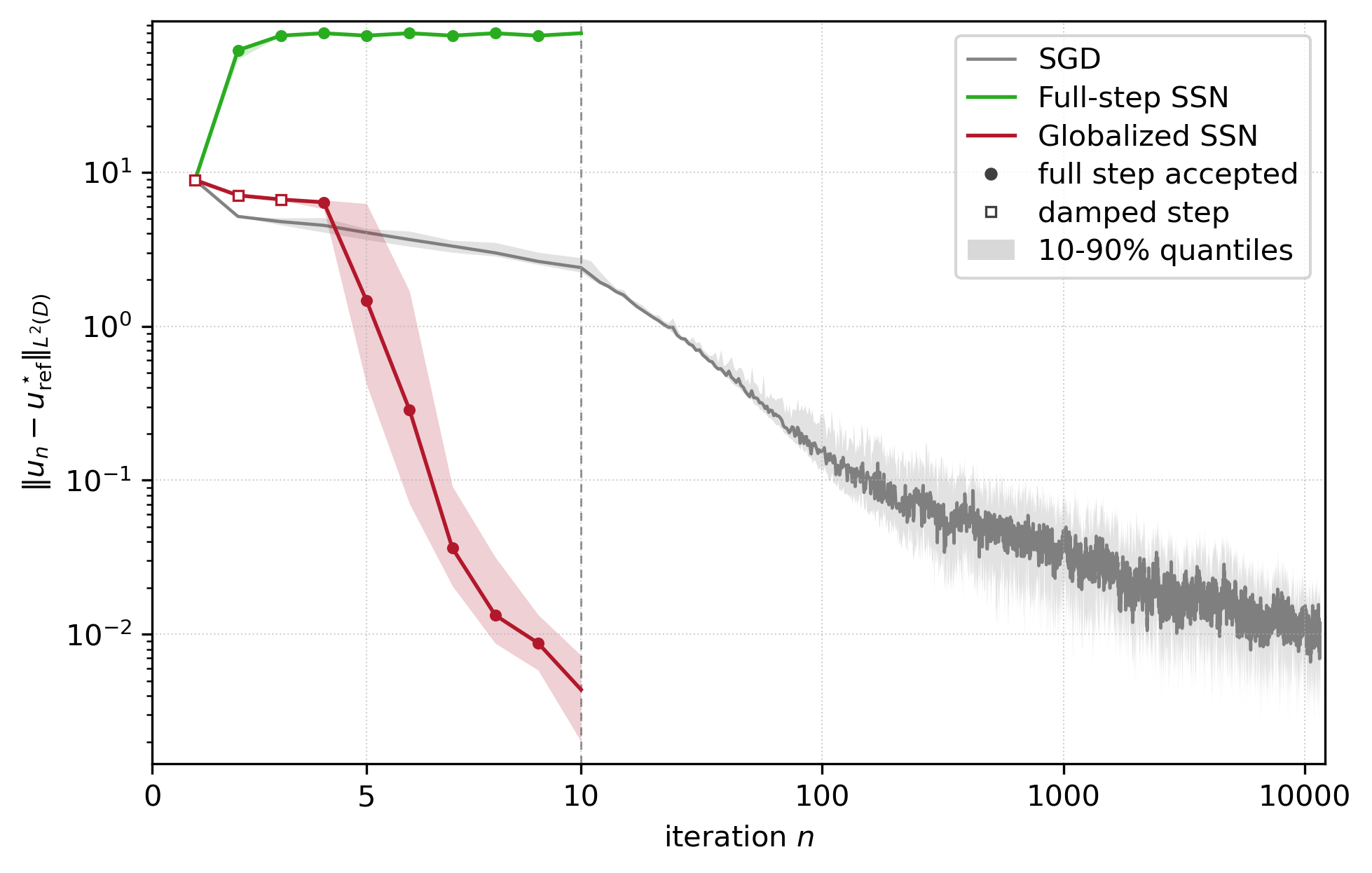}
        \caption{\centering Convergence vs. iterations}
        \label{fig:SSSN_elliptic_iterations}
    \end{subfigure}
    \hfill
    \begin{subfigure}{0.48\textwidth}
        \centering
        \includegraphics[width=\linewidth]{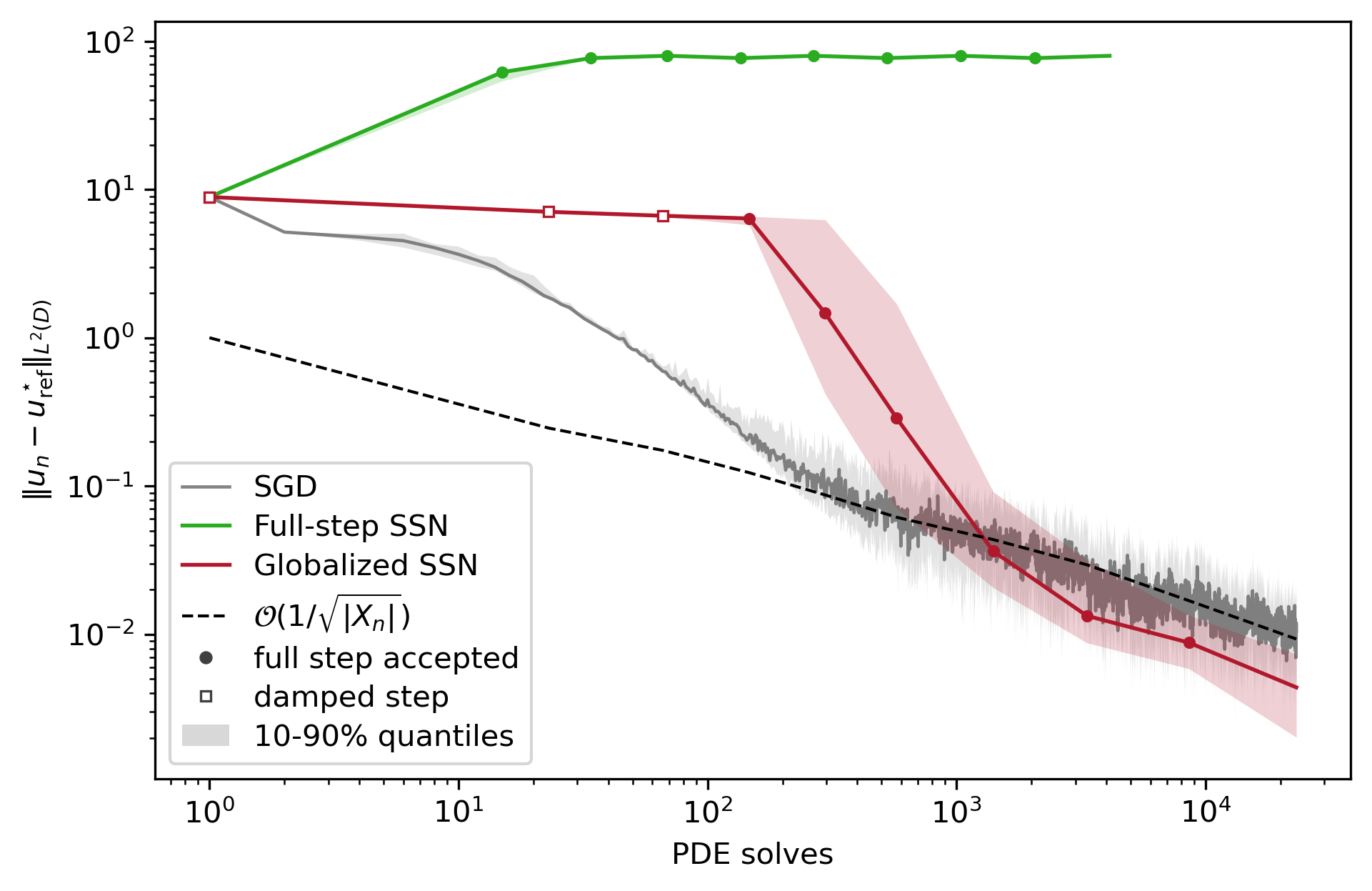}
        \caption{\centering Convergence vs.~PDE solves}
        \label{fig:SSSN_elliptic_PDE_solves}
    \end{subfigure}
    \caption{Convergence of iterates in the $L^2(D)$-norm for SGD, full-step SSN, and globalized SSN for the elliptic problem.}
    \label{fig:comparison_elliptic}
\end{figure}

We run each method a total of 20 times and track the averaged errors, i.e., $\frac{1}{20}\sum_{i=1}^{20}\|u^i_n - u^\ast_{\text{ref}}\|_{L^2(D)} \approx \EE{\|u_n - u^\ast_{\text{ref}}\|_{L^2(D)}} $
to account for the inherent randomness of the problem. We also track the $10-90\%$ quantiles of the error, i.e., the band in which $80\%$ of the errors lie. We also track the SSN-gradient induced noise floor $\mathcal{O}(1/\sqrt{|X_n|}),$ as this can also be seen as a bound for the rate of convergence of \cref{alg:globalized_SSN}.
The results are displayed in \cref{fig:comparison_elliptic}. We note that in \cref{fig:SSSN_elliptic_iterations}, the scale for iteration $n$ is linear on $[0, 10]$ and logarithmic beyond to better differentiate between SSN and SGD.

Comparing full-step SSN with globalized SSN in \cref{fig:comparison_elliptic}, we observe that the full-step SSN first moves away from the solution in its first iterate and then stays at this plateau.
The globalized SSN, however, takes three damped steps on average before switching to full Newton steps and converging to the Monte Carlo floor. In the local phase, the error of \cref{alg:globalized_SSN} decreases to the level of the noise floor $\mathcal{O}(1/\sqrt{|X_n|})$
in \cref{fig:SSSN_elliptic_PDE_solves}. Thus, the attainable accuracy is also limited by the gradient batch size.
We can observe that the switch to the local phase happens once we are predominantly in the quadratic region of $\Psi_r$, i.e., once the $L^2(D)$-tracking is active on most of $D$. The ratio of $x\in D$ where the $L^2(D)$ tracking is active can be estimated by 
\begin{equation}\label{eq:quadratic_ratio}
      \frac{1}{\operatorname{meas}(D)}
    \int_D \mathbb{P}\big(|y(x,\xi) - y_d(x)| \leq r\big)
    \,\mathrm{d}x.
\end{equation}
In \cref{tab:huber_quadratic_fraction} we report the estimation of \eqref{eq:quadratic_ratio} at iterates $u_n$ as well as the ratio of accepted tests \eqref{eq:ast}, i.e., full steps at $u_n$ of the 20 runs. We see that the $L^2(D)$ tracking fraction grows from $0.31$ to above $0.9$ during the damped phase and that the full step is accepted once this fraction exceeds roughly $0.9$ (iterations $n = 3$ to $5$). From $n=5$ onwards, all runs accept the full step and the Huber term is quadratic on essentially all of $D$. Hence, near $u^\ast$, the problem reduces to a linear--quadratic one and the example mainly illustrates the need for globalization in the initial phase. In the Moreau--Yosida example below, the nonsmooth term remains active along the entire iteration sequence, specifically at the reference solution (see \cref{tab:gamma_sweep_parabolic}). 
\begin{table}[htbp]
  \centering
  \caption{Estimation of \eqref{eq:quadratic_ratio} at the iterates $u_n$
    of \cref{alg:globalized_SSN},
    and fraction of the $20$ runs
    accepting the full step at $u_n$.}
  \label{tab:huber_quadratic_fraction}
  \small
  \setlength{\tabcolsep}{5pt}
  \begin{tabular}{@{}l*{10}{c}@{}}
    \toprule
    iteration $n$
      & 0 & 1 & 2 & 3 & 4 & 5 & 6 & 7 & 8 & 9 \\
    \midrule
    $L^2(D)$ tracking fraction
      & 0.31 & 0.51 & 0.72 & 0.89 & 0.94
      & 1.00 & 1.00 & 1.00 & 1.00 & 1.00 \\
    full step fraction
      & 0.00 & 0.00 & 0.10 & 0.55 & 0.95
      & 1.00 & 1.00 & 1.00 & 1.00 & -- \\
    \bottomrule
  \end{tabular}
\end{table}

Comparing \cref{alg:globalized_SSN} to SGD, we note that while \cref{alg:globalized_SSN} takes only 9 iterations and SGD takes 11577 iterations, the two methods reach a similar error in terms of PDE solve complexity. This matches \cref{rem:oracle_complexity}. 
With geometrically growing batches, the global phase requires $\mathcal{O}(\varepsilon^{-1})$ PDE solves, the same order as SGD, but only $\mathcal{O}(\log(\varepsilon^{-1}))$ iterations.
The advantage of \cref{alg:globalized_SSN} therefore lies in the number of iterations, whose samples can be processed in parallel, rather than in the total work.
In numerical experiments, using mini-batched SGD with a fixed or growing batch-size always seemed to perform worse than single sample SGD with respect to PDE solve complexity.

\subsection{Linear parabolic PDE with Moreau--Yosida regularization}
Next, we study an optimal control problem subject to a parabolic PDE on the space-time cylinder $Q := D \times (0,T)$ for a given $T>0$. We consider almost sure state constraints of type $y(\cdot, \xi) \leq y_c$ a.e.~in $Q$, a.s.~in $\Xi$. Since the associated Lagrange multipliers are measures \cite[p.~340]{troltzsch2010optimal}, these state constraints are notoriously difficult to handle and are frequently regularized by a Moreau--Yosida penalty term 
$\frac{\gamma}{2}\mathbb{E}[{\|\max(0, y(\cdot, \xi)-y_c)\|^2_{L^2(Q)}}]$ with $\gamma>0$ and $L^2(Q):= L^2(0, T; L^2(D)),$  see e.g.~\cite{geiersbach2022OptimalityConditionsMoreau, hintermuller2009MoreauYosidaRegularization, hintermuller2006FeasibleNoninteriorPathFollowing, hintermuller2006pathfollowing}, and for related augmented Lagrangian approaches see e.g.~\cite{ito2000AugmentedLagrangianNonsmooth,ito2003SemismoothNewtonMethods}.

\subsubsection{Setup}
We now consider $D = (0,1),$ $T = 1$ and the optimization problem 
\[\min_{u \in L^2(Q)} j(u) :=  \frac{1}{2}\EE{ \|y(\cdot, \xi)-y_d\|^2_{L^2(Q)} + {\gamma} \|\max(0, y(\cdot, \xi)-y_c)\|^2_{L^2(Q)}} + \frac{\mu}{2}\|u\|_{L^2(Q)}^2\]
subject to the random linear parabolic state equation 
\begin{equation*}\label{eq:parabolic_example_state}
    \begin{aligned}
        \partial_t y(x,t, \xi) - \partial_x ( a(x, \xi)\, \partial_x y(x,t,\xi) ) &= u(x,t) + f(x, t, \xi), && (x,t) \in Q,\ \xi \in \Xi, \\
y(0, t, \xi) = y(1, t, \xi) &= 0, && t \in (0,T),\\
y(x, 0, \xi) = y_0(x) &= 0, && x \in D.\\
    \end{aligned}
\end{equation*}
For the (time-independent) diffusion coefficient, we use the model of \cref{sec:linear_elliptic} with $\psi_k(x) = \cos(k \pi x)$, $k=1,\dots,4$, which again gives $a(x, \xi) \in [0.45, 1.55].$
For the right-hand side, we set \[f(x,t,\xi) = 0.5\,\big(\xi^r_1\, \sin(\pi x) + \xi^r_2\, \sin(2\pi x)\sin(\pi t / T)\big)\]
 with i.i.d.~$\xi^r_k  \, \sim \,  \mathcal{U}(-1,1),$ which again gives $ \xi = (\xi_1^l, \xi_2^l, \xi_3^l, \xi_4^l, \xi_1^r, \xi_2^r) \in \Xi \subset \mathbb{R}^6$.
Further, we set \begin{align*}
    y_d(x,t) = 2 \sin( \pi t /T)\sin(\pi x), \qquad y_c(x,t) \equiv 0.1.
\end{align*}

For first-order optimality, the adjoint equation runs backward in time and is given by 
\begin{equation*}
  \begin{aligned}
    -\partial_t p - \partial_x\bigl(a\,\partial_x p\bigr)
      &= y - y_d + \gamma\,\max(0, y - y_c)
      && \text{in } Q \times \Xi, \\
    p(0,\cdot,\cdot) = p(1,\cdot,\cdot) &= 0
      && \text{on } (0,T) \times \Xi, \\
    p(\cdot,T,\cdot) &= 0
      && \text{in } D \times \Xi,
  \end{aligned}
\end{equation*}
with the Riesz gradient in $L^2(Q)$ given by 
\[g(u) = \nabla j(u) = \mathbb{E}[p(\cdot, \cdot, \xi)] + \mu u,\]
(see e.g.~\cite{troltzsch2010optimal}). 
Verifying \cref{ass:lipschitz_mu_convex,ass:newton_differentiability,ass:gradient_variance,ass:hessian_variance,ass:hessian_bounds} is mostly analogous to the elliptic case in \cref{sec:linear_elliptic}. We note that now, for each $\xi \in \Xi$ the corresponding state $y(\cdot, \xi)$ belongs to the space \[W(0,T):= \{y : y \in L^2(0,T; H^1_0(D)), y_t \in L^2(0,T; H^{-1}(D))\}\hookrightarrow L^2(Q).\] This motivates the bounded operator $S(\xi): L^2(Q) \to L^2(Q)$ with 
\[\|S(\xi)\|_{L^2(Q) \to L^2(Q)} \leq \bar{S} \quad \text{for all }\xi \in \Xi.\]
Further, we know that $W(0,T) \hookrightarrow C([0,T]; L^2(D))$
\cite[Theorem 3.10]{troltzsch2010optimal}, and hence $W(0,T) \hookrightarrow L^\infty(0,T; L^2(D)) \cap L^2(0,T; H^1_0(D))$. By \cite\cite[Chapter I, Proposition 3.1]{dibenedetto1993DegenerateParabolicEquations}, we obtain $W(0,T) \hookrightarrow L^q(Q)$ with $q = 2(d+2)/d > 2$ for every $d \geq 1$, which closes the norm gap required for Newton differentiability of the $\max$ operator.
Lastly, we note that in this example, we have affinely growing gradient variance with $\nu_1 \neq 0$ as stated in \cref{ass:gradient_variance}. 
Unlike the elliptic example where $\min(\max(y(x, \xi) - y_d, -r), r)$ could be uniformly bounded from above and below, this does not hold for a right hand side $\max(0, y(x, t, \xi) - y_c(x,t))$. Hence, with \eqref{eq:variance_inequality}, computing the gradient variance results in 
\begin{align*}
    \mathbb{E}[\|G_{\xi}(u) - g(u)\|_{L^2(Q)}^2] &=\mathbb{E}[\|p(\cdot, \xi)  - \mathbb{E}[p(\cdot, \xi)] \|_{L^2(Q)}^2] \leq \mathbb{E}[\|p(\cdot, \xi)  \|_{L^2(Q)}^2] \\
     &\leq \bar{S}^2\mathbb{E}[\|y(\cdot, \xi) - y_d + \gamma\max(0, y(\cdot, \xi) - y_c) \|_{L^2(Q)}^2]\\
     &\leq \bar{S}^2\mathbb{E}\left[ \left( (1 + \gamma) \bar{S} (\|u\|_{L^2(Q)} + \|f(\cdot, \xi)\|_{L^2(Q)}) \right.\right.\\
     & \quad \left. \left.+ \|y_d\|_{L^2(Q)} + \gamma \|y_c\|_{L^2(Q)} \right)^2\right]\\
    & \leq \nu_0^2 + \nu_1^2\|u- u^\ast\|_{L^2(Q)}^2
\end{align*}
for appropriately chosen $\nu_0^2, \nu_1^2$ which, we note, are of order $\gamma^2$. 

\paragraph{Path following} In practice, the Moreau--Yosida parameter $\gamma$ is driven to infinity along a sequence $\gamma_k \to \infty$. In the deterministic setting, for each $\gamma_k$, the resulting problem is usually solved with a semismooth Newton method and warm-started with the approximate previous solution for $\gamma_{k-1}$, and $\gamma_k$ is updated based on a value functional \cite{hintermuller2006FeasibleNoninteriorPathFollowing, hintermuller2006pathfollowing}. 
In the stochastic setting with almost sure constraints, consistency of the regularization as $\gamma \to \infty$ was explored in \cite{geiersbach2022OptimalityConditionsMoreau}. In the following example, we will use a heuristic version of this path following. We fix $\gamma_k = 10^k$, $k = 1,2,3,$ and for each $\gamma_k$ re-estimate $L$ and $M$ as in \cref{sec:linear_elliptic} and hence the test parameter $\theta$ and the damped step $m/L$. Using warm-starts from the previous solution for $\gamma_2, \gamma_3$, we skip the global phase and resume the local batch schedule. Yet, we will use the same local batch schedule for each $\gamma_k$. A stochastic path-following method that couples the update of $\gamma_k$ to the batch sizes is left for future work.

\subsubsection{Numerical results}
We discretize the space with piecewise linear, globally continuous finite elements on a uniform grid of $D$ with $10$ cells and in time with the implicit Euler method with $10$ steps, resulting in $90$ degrees of freedom. The Newton systems $H_{S_n}(u_n)[d_n] = -G_{X_n}(u_n)$ are again solved with CG in the $L^2(Q)$ inner product with relative tolerance $10^{-2}$.

We solve the problem for $\gamma \in \{10, 10^2, 10^3\}$. For $\gamma = 10$, we start from $u_0 \equiv 50$ and each subsequent problem is warm-started from the solution of the previous $\gamma$. Throughout, we set $\mu = m = 10^{-1}$ and re-evaluate $M = L$ as in \cref{sec:linear_elliptic}. For each $\gamma,$ a SAA reference solution $u^\ast_{\text{ref}, \gamma}$ is computed on an independent sample set of size $N = 10^5$ with a deterministic semismooth Newton method. This implies a reference accuracy $\|u^\ast_{\gamma} - u^\ast_{\text{ref}, \gamma}\|_{L^2(Q)}\approx 1/\sqrt{N} \approx 3\cdot 10^{-3}.$ This limited accuracy along with the sampling error makes it difficult to observe superlinear convergence in the error plots. 

We use the batch schedule of \cref{tab:schedules_elliptic} with $q = \sqrt{2/3}$. For $\gamma = 10$, we switch to the local phase as in \cref{sec:linear_elliptic} and run six local iterations. 
The warm-started problems for $\gamma = 10^2, 10^3$ skip the global phase and resume the local schedule at $k=3$ for three more iterations, since restarting from small batches would place the sampling error of the first Newton step far above the accuracy of the warm start. 
Again, as in \cref{sec:linear_elliptic} we compare with single sample SGD run with step size $\tau_n = \frac{1}{(n+1)\mu}$ for the same amount of PDE solves. We run both methods $20$ times per $\gamma$ and report the mean and the $10$--$90\%$ quantiles of
$\|u_n -u^\ast_{\text{ref},\gamma}\|_{L^2(Q)}$.

\begin{figure}[ht]
  \centering
  \includegraphics[width=\linewidth]{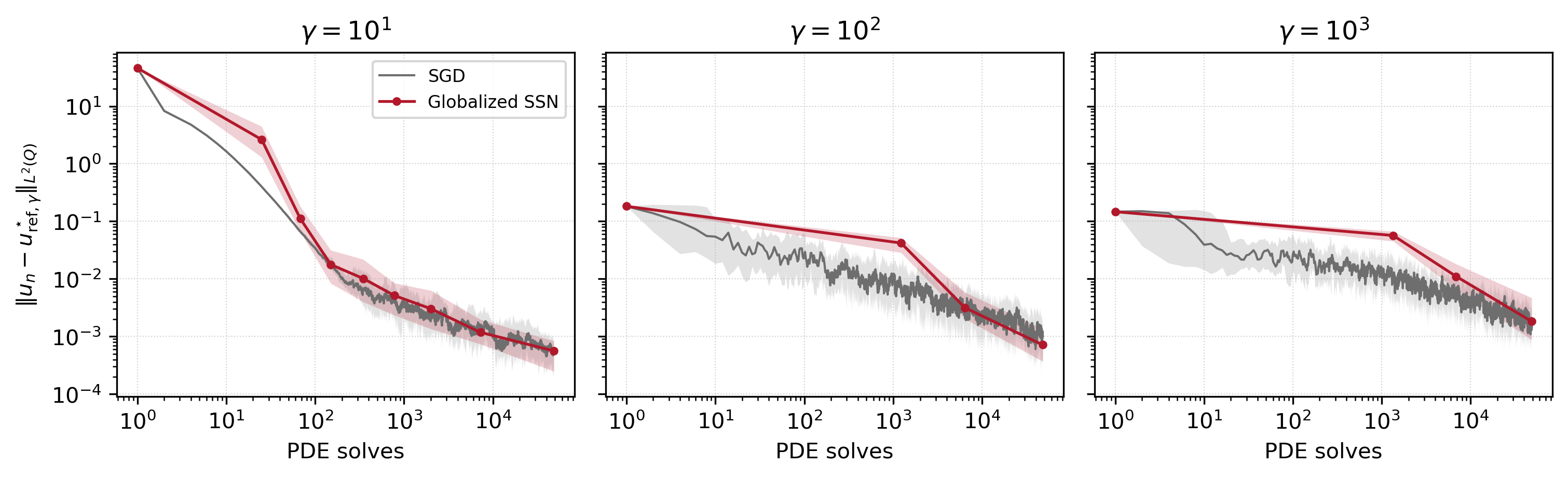}
  \caption{Average and $10$--$90\%$ quantiles
    of $\|u_n - u^\ast_{\text{ref},\gamma}\|_{L^2(Q)}$
    versus PDE solves over $20$ runs,
    for $\gamma = 10 \to 10^2 \to 10^3$.}
  \label{fig:gamma_sweep_parabolic}
\end{figure}

\begin{table}[htbp]
  \centering
  \caption{Average number of damped SSN steps per run, total iterates of SSN, PDE solves, violation $\mathbb{E}[\|\max(0, y_\gamma(\cdot, \xi) - y_c) \|_{L^2(Q)}]$ and active fraction 
  $\frac{1}{\operatorname{meas}(Q)}\int_Q \mathbb{P}(y(x,t,\xi) > y_c(x,t))\,\mathrm{d}(x,t)$ for different values of $\gamma$.}
  \label{tab:gamma_sweep_parabolic}
\begin{tabular}{cccccc}
  \toprule
  $\gamma$ & damped steps &  iter. SSN & PDE solves &  violation & active fraction  \\
  \midrule
  $10$ & 0.1 & 8 & $4.8\cdot 10^{4}$ &  $1.2\cdot 10^{-2}$  & 0.24\\
  $10^{2}$  & 0.0  & 3 & $4.8\cdot 10^{4}$ &  $3.9\cdot 10^{-3}$ & 0.12\\
  $10^{3}$ & 0.0 &  3 & $4.9\cdot 10^{4}$  & $6.6\cdot 10^{-4}$ & 0.04\\
  \bottomrule
\end{tabular}
\end{table}

\Cref{fig:gamma_sweep_parabolic} and \cref{tab:gamma_sweep_parabolic}
show the following effects. First, similar to \cref{sec:linear_elliptic}, \cref{alg:globalized_SSN} takes a similar amount of PDE solves compared to SGD to reach the same accuracy. Note that \cref{alg:globalized_SSN} only produces an iterate at the end of a Newton step, so the long segments between markers in \cref{fig:gamma_sweep_parabolic} reflect the cost of one large batch, not stagnation.
Second, damped steps are rare: only two damped steps occur over the $20$ runs, both in the cold start $u_0 \equiv 50$ for $\gamma = 10$. The warm-started problems for $\gamma = 10^2, 10^3$ start close to $u^\ast_{\mathrm{ref},\gamma}$ and accept every full step.
Third, the expected violation of the state constraint decreases with $\gamma$,
consistent with the Moreau--Yosida consistency results in \cite{geiersbach2022OptimalityConditionsMoreau}. 
Fourth, the active fraction $\frac{1}{\operatorname{meas}(Q)}\int_Q \mathbb{P}(y(x,t,\xi) > y_c(x,t))\,\mathrm{d}(x,t)$ computed at the reference solution $u^\ast_{\mathrm{ref}, \gamma}$ decreases from $0.24$ to $0.04$, but remains positive. In contrast to \cref{sec:linear_elliptic}, the nonsmooth term is therefore active at the solution. 
As in the elliptic example, all samples of one Newton iteration can be evaluated in parallel whereas SGD is strictly sequential.

\section{Conclusion}
In this work, we proposed a stochastic semismooth Newton method in Hilbert spaces for strongly convex minimization problems with mini-batched gradients and generalized Hessians and affinely growing stochastic gradient variance. 
We proved global R-linear convergence in expectation which, almost surely, transitions to R-superlinear convergence with full Newton steps, without conditions on the initial point or generated sequence. 
Rates and batch sizes depend only on the problem constants of the infinite-dimensional problem. Numerical experiments showed that the method succeeds, where the full-step method can fail. In terms of PDE-solve cost, the method matched that of SGD with only a few iterations that can be evaluated in parallel. Open directions include a stochastic path following method that couples the batch sizes to the Moreau--Yosida parameter, and studying control constraints within a stochastic augmented Lagrangian method. 

\bibliographystyle{siamplain}
\bibliography{references}

\end{document}